%% file: main_arxiv.tex
\documentclass[journal]{IEEEtran}
\usepackage{amsmath, amssymb, epsfig, graphicx, bm}
\usepackage{cases}
\usepackage{float, subfigure}
\usepackage{color}
\usepackage{amsfonts}
\usepackage{algorithm}
\usepackage{algpseudocode}
\usepackage{cite, url}
\usepackage[all,cmtip]{xy}
\usepackage{color,hyperref}
\definecolor{darkblue}{rgb}{0,0,0.8}
\hypersetup{colorlinks,breaklinks,
	linkcolor=darkblue,urlcolor=darkblue,anchorcolor=darkblue,citecolor=darkblue}

\newtheorem{theorem}{Theorem}
\newtheorem{lemma}{Lemma}
\newtheorem{assumption}{Assumption}

\newtheorem{remark}{Remark}

\input{notation}
\title{\LARGE \bf 
	A Compressed Gradient Tracking Method for Decentralized Optimization with Linear Convergence
}

\author{Yiwei Liao, Zhuorui Li,  Kun Huang, and Shi Pu% <-this % stops a space
	\thanks{This work was partially supported  by the Shenzhen
		Research Institute of Big Data (SRIBD) (Grant No. J00120190011), by the National Natural Science Foundation of China (NSFC) (Grant No. 62003287), by the Shenzhen Science and Technology Program (Grant  No. RCYX20210609103229031 and No. GXWD20201231105722002-20200901175001001), and by Shenzhen Institute of Artificial Intelligence and Robotics for Society (AIRS) (Grant No. AC01202101108). Yiwei Liao  and Zhuorui Li contributed equally to this work. 
		Corresponding author: Shi Pu.}% <-this % stops a space
	\thanks{Yiwei Liao is with the School
		of Data Science, The Chinese
		University of Hong Kong, Shenzhen, China and also with the Shcool of Information Science and Technology, University of Science and Technology of China, Hefei, China. Zhuorui Li is with Shenzhen Research Institute of Big Data, Shenzhen, China. Kun Huang is with the School
		of Data Science, Shenzhen Research Institute of Big Data, The Chinese
		University of Hong Kong, Shenzhen, China. Shi Pu is with the School of Data Science, Shenzhen Institute of Artificial Intelligence and Robotics for Society (AIRS), The Chinese
		University of Hong Kong, Shenzhen, China.
		{\tt\small (emails: lyw@stu.scu.edu.cn, lizhuorui27@gmail.com, kunhuang@link.cuhk.edu.cn, pushi@cuhk.edu.cn)}}%
}

\begin{document}
	\maketitle
	\thispagestyle{empty}
	\pagestyle{empty}
	
	%%%%%%%%%%%%%%%%%%%%%%%%%%%%%%%%%%%%%%%%%%%%%%%%%%%%%%%%%%%%%%%%%%%%%%%%%%%%%%%%%%%%%%%%%%%%%%%%%%	
	\begin{abstract}
		Communication compression techniques are of growing interests for solving the decentralized optimization problem under limited communication, where the global objective is to minimize the average of local cost functions over a multi-agent network using only local computation and peer-to-peer communication. In this paper, we   propose a novel compressed gradient tracking algorithm (C-GT) that combines gradient tracking technique with communication compression. In particular, C-GT is compatible with a general class of compression operators that unifies both unbiased and biased compressors. We show that C-GT inherits the advantages of gradient tracking-based algorithms and achieves linear convergence rate for strongly convex and smooth objective functions.  
		Numerical examples complement the theoretical findings and demonstrate the efficiency and flexibility of the proposed algorithm.
	\end{abstract}
	

	%%%%%%%%%%%%%%%%%%%%%%%%%%%%%%%%%%%%%%%%%%%%%%%%%%%%%%%%%%%%%%%%%%%%%%%%%%%%%%%%%%%%%%%%%%%%%%%%%%
	\section{Introduction}
	In this paper, we study the problem of decentralized optimization over a multi-agent network that consists of  $n$ agents. The goal is to collaboratively solve the following optimization problem:
	\begin{equation} \label{problem}
		\begin{array}{c}
			\min\limits_{x\in\RR^p}~f(x):=\frac{1}{n}\sum\limits_{i=1}^n f_i(x),
		\end{array}
	\end{equation}
	where $x$ is the global decision variable, and each agent has a local objective function $f_i: \RR^p\rightarrow \RR$. The agents are connected through a communication network and can only exchange information with their immediate neighbors in the network. Through local computation and local information exchange, they seek a consensual and optimal solution that minimizes the average of all the local cost functions. 
	Decentralized optimization is widely applicable when central controllers or servers are not available or preferable, when centralized communication that involves a large amount of data exchange is prohibitively expensive due to limited communication resources, and when privacy preservation is desirable.
	
	Problem \eqref{problem} has  attracted much attention in recent years and has
	found a variety of  applications in wireless networks, distributed control of robotic systems, and machine learning, etc  \cite{Cohen2016Distributed,Nedic2018Distributed,Nedic2020Distributed}. 
	To solve \eqref{problem} over a multi-agent network, early work considered the distributed subgradient descent (DGD) method with a diminishing step-size policy \cite{Nedic2009distributed}. 
	Under a constant step-size, EXTRA \cite{Shi2015Extra} first achieved linear convergence rate for strongly convex and smooth cost functions by introducing an extra correction term to DGD. Distributed gradient tracking-based methods were later developed in \cite{Xu2015Augmented,Di2016Next,Nedic2017achieving,Qu2018Harnessing}, where the local gradient descent direction in DGD was replaced by an auxiliary variable that is able to track the average gradient of the local objective functions.  As a result, each agent's local iteration  is moving in the global descent direction and converges exponentially  to the optimal
	solution for strongly
	convex and smooth objective functions \cite{Nedic2017achieving,Qu2018Harnessing}.  Compared with EXTRA, gradient tracking-based methods are also suitable for uncoordinated step-sizes \cite{Xu2015Augmented,Nedic2017Geometrically} and possibly asymmetric weight matrices while preserving linear convergence rates.  
	Some variants were also proposed to deal with stochastic gradient information and time-varying or directed network topologies, etc. For example, in \cite{Pu2020distributed}, a distributed stochastic gradient tracking method was considered which exhibits comparable performance to a centralized stochastic gradient algorithm.  Combining an approximate Newton-type method and gradient tracking leads to Network-DANE, which enables further computational savings by performing variance-reduced techniques \cite{Li2020communication}. 
	Time-varying networks were considered in \cite{Nedic2015Distributed,Nedic2017achieving,Xie2018Distributed,sun2022distributed}, and  
	more recent development on directed graphs can be found in \cite{Tsianos2012PushSum,Nedic2015Distributed,Xin2018Linear,Pu2021Push,Xin2020General,Pu2020Robust} and the references therein.

	In many application scenarios,  it is vital to design  communication-efficient protocols for distributed computation due to limited communication bandwidth and power constraints. Recently, in order to improve system scalability and communication efficiency, researchers have considered a variety of communication compression methods, such as sparsification and quantization \cite{1-bit-sgd, NIPS2017_qsgd, signSGD-ICML18, Stich2018Sparsified, karimireddy2019Error, mishchenko2019distributed, tang2019doublesqueeze, Stich2020Communication,Beznosikov2020Biased,Xu2020Compressed},   under the master-worker centralized architecture.  
	Several techniques were  introduced to alleviate compression errors, including   compression error compensation and gradient difference compression \cite{1-bit-sgd,Stich2018Sparsified,mishchenko2019distributed,tang2019doublesqueeze}. 
	
	In the decentralized setting, the difference compression and extrapolation compression techniques were introduced to reduce model compression error in \cite{tang_NIPS2018_7992}.  A novel algorithm with communication compression (CHOCO-SGD), which combines with DGD and preserves the model average,  was presented in \cite{ksj2019choco, koloskova*2020decentralized}. But the method converges sublinearly even when the objective functions are strongly convex. In \cite{liu2020linear}, a linearly convergent decentralized optimization algorithm  with compression (LEAD) was introduced for strongly convex and smooth objective functions. The method is based on NIDS \cite{li2019decentralized}, a sibling of EXTRA. In light of an incremental primal-dual method, a linearly convergent quantized decentralized optimization algorithm was developed for   unbiased randomized compressors in \cite{kovalev2021linearly}. In \cite{Lee2021finite}, a black-box model was provided for distributed algorithms
	based on finite-bit quantizer. 
	
	In light of the advantages of gradient tracking-based methods for decentralized optimization, it is natural to consider the marriage between gradient tracking and communication compression. The first such effort was made in \cite{Kajiyama2020Linear} which considered a quantized gradient tracking method based on a special quantizer. It was shown to achieve linear convergence rate for strongly convex and smooth objective functions. However, the algorithm design is rather complicated and relies on a specific quantizer. In addition, the convergence conditions are not easy to verify. 
	
	In this paper, we   consider  a novel gradient tracking-based method (C-GT) for decentralized optimization with communication compression. The algorithm compresses both the decision variables and the gradient trackers to provide a communication-efficient implementation. Unlike the existing methods which are mostly based on unbiased compressors or biased but contractive compressors, C-GT is provably efficient for a general class of compressors, including those which are neither unbiased nor biased but contractive,  e.g., the composition of quantization and sparsification and the norm-sign compression operators. We show that C-GT achieves linear convergence for strongly convex and smooth objective functions under such a general class of communication compression techniques,  where the agents may choose different, uncoordinated step-sizes.
	
	The main contributions of the paper are summarized as follows:
	\begin{itemize}
		\item We propose a novel compressed gradient tracking  algorithm (C-GT)  for decentralized optimization, which inherits the advantages of gradient tracking-based methods and saves communication costs at the same time.
		\item The proposed C-GT algorithm is applicable to a general class of compression operators and  works under arbitrary compression precision. In particular, the general condition on the compression operators unifies the commonly considered unbiased and biased but contractive compressors and also includes other compression methods such as  the composition of quantization and sparsification and the norm-sign compressors.
		\item C-GT provably achieves linear convergence for minimizing strongly convex and smooth objective functions under the general condition on the compression operators,  where the agents may choose different, uncoordinated step-sizes.
		\item Simulation examples show that C-GT is efficient  compared to the state-of-the-art methods and widely applies to various compressors.
	\end{itemize}
	
	The rest of this paper is organized as follows. 
	We present the general condition on the compressors and the C-GT algorithm  in Section \ref{sec: Alg}. In Section \ref{sec: CA-CGT}, we perform the convergence analysis for C-GT. 
	Numerical examples are provided in Section \ref{sec: simulation}. Finally, concluding remarks are given in Section \ref{sec: conclusion}.
	\subsection{Notation}\label{subsec: Notation}
	Vectors are columns if not otherwise specified in this paper.
	Let each agent $i$ hold a local copy $\vx_i\in\mathbb{R}^p$ of the decision variable and a gradient tracker (auxiliary variable) $\vy_i\in\mathbb{R}^p$. At the $k$-th iteration, their values are denoted by $\vx_{i}^{k}$ and $\vy_{i}^{k}$, respectively. 
	For notational convenience, define
	$\vX := [\vx_1, \vx_2, \ldots, \vx_n]^{\T}\in\mathbb{R}^{n\times p}$, 
	$ \vY := [\vy_1, \vy_2, \ldots, \vy_n]^{\T}\in\mathbb{R}^{n\times p}$,
	and
	$\oX :=  \frac{1}{n}\mathbf{1}^{\T} \vX\in\mathbb{R}^{1\times p}$,  
	$\oY :=  \frac{1}{n}\mathbf{1}^{\T}\vY\in\mathbb{R}^{1\times p}$, 
	where $\mathbf{1}$ is the column vector with each entry given by 1. At the $k$-th iteration, their values are denoted by $\vX^{k}$, $\vY^{k}$, $\oX^{k}$ and $\oY^{k}$, respectively. 
	Auxiliary variables of the agents (in an aggregative matrix form) $\vH_{x}$, $\vH_{y}$, $\vQ_{x}$, $\vQ_{y}$, $\widehat{\vX}$,and $\widehat{\vY}$ are defined similarly. 
	Denote the aggregative gradient 
	$
	\nabla \vF(\vX):=\left[\nabla f_1(\vx_1), \nabla f_2(\vx_2), \ldots, \nabla f_n(\vx_n)\right]^{\T}\in\mathbb{R}^{n\times p},
	$ 
	and 
	$
	{\nabla} \overline {\vF}(\vX) := \frac{1}{n}{\vone^\T \nabla \vF(\vX)}=\frac{1}{n}\sum_{i=1}^n\nabla f_i(\vx_i).
	$
	
	We use $\|\cdot\|$ to denote the Frobenius norm of vectors and matrices by default. Specially, for square matrices, $\|\cdot\|$ represents the spectral norm. The spectral radius of a square matrix $\vM$ is denoted by $\rho(\vM)$.
	
	\section{Problem formulation}\label{sec: Pre}
	In this section, we provide the assumptions on the communication graphs and the objective functions.  Then, we discuss different kinds of  compression methods and provide a general description for compression operators.
	
	\subsection{Preliminaries}
	We start with introducing the conditions on the communication network/graph and the objective functions. 
	Assume the agents are connected over an undirected graph   $\mathcal{G}=(\mathcal{V},\mathcal{E})$, where $\mathcal{V}=\{1,2,\ldots,n\}$ is the set of vertices (nodes) and $\mathcal{E}\subseteq \mathcal{V}\times \mathcal{V}$  is the set of edges.  For an arbitrary agent $i\in\mathcal{V}$, we define the set of its neighbors as $\mathcal{N}_{i}$.  
	Regarding the network structure, we make the following standing assumption.
	\begin{assumption}\label{Assumption: network}
		The undirected graph $\mathcal{G}$ is  strongly connected and permits a nonnegative doubly stochastic weight matrix $\vW=[w_{ij}]\in\RR^{n\times n}$. That is, agent  $i$ can receive  information from agent $j$ if and only if $w_{ij}>0$, and $\vW \vone=\vone$  and $\vone ^{\T}\vW=\vone ^{\T}$.
	\end{assumption}

	\begin{remark}
		Although we assume an undirected graph $\mathcal{G}$, note that the considered C-GT method also works with any balanced directed graph, where it is convenient to construct a doubly stochastic weight matrix.
	\end{remark}

	The assumption on the objective functions is given below.
	
	\begin{assumption}
		\label{Assumption: function}
		The local cost function $f_i$ is $\mu_i$-strongly convex, and 
		its gradient is $L_i$-Lipschitz continuous, i.e., for any $\vx,\vx'\in\mathbb{R}^p$,
		\begin{align}
			& \langle \nabla f_i(\vx)-\nabla f_i(\vx'),\vx-\vx'\rangle\ge \mu_i\|\vx-\vx'\|^2,\\
			& \|\nabla f_i(\vx)-\nabla f_i(\vx')\|\le L_i \|\vx-\vx'\|.
		\end{align}
	\end{assumption}
	From Assumption \ref{Assumption: function}, the objective function $f$ is $\mu$-strongly convex and the gradient of $f$ is $L$-Lipschitz continuous, where $\mu=\frac{1}{n}\sum_{i=1}^n \mu_i$ and $L=\max{\{L_i\}}$.\footnote{We denote $\kappa=L/\mu$ as the condition number.}
	Moreover, there exists a unique solution denoted by $\vx^*\in\mathbb{R}^{1\times p}$  to problem \eqref{problem}   under Assumption \ref{Assumption: function}.

	\subsection{Compression Methods}
	In this subsection, we introduce some common assumptions on the compression operators and then present a more general and unified assumption.
	\subsubsection{Unbiased compression operators}
	\begin{assumption}\label{Assumption:Unbiased}
		The  compression operator $\cQ:\RR^d\rightarrow\RR^d$  satisfies $\EE \cQ(\vx)=\vx$ and there exists a constant $C\geq0$ such that $\EE\|\cQ(\vx)-\vx \|^2\leq C\|\vx\|^2$, $\forall \vx \in \RR^d$.
	\end{assumption}
	
	\begin{remark}
		The expectation  is taken with respect to the random vector corresponding to the internal compression randomness of $\cQ$. 
		Some instances of feasible stochastic
		compression operators satisfying Assumption \ref{Assumption:Unbiased}, such as the  unbiased $b$-bits $q$-norm quantization compression method, can be found in \cite{ksj2019choco, koloskova*2020decentralized,liu2020linear}  and the references therein. 
	\end{remark}
	\subsubsection{Biased compression operators}
	\begin{assumption}\label{Assumption:Biased}
		The  compression operator $\cC_\delta \colon \RR^d \to \RR^d$   satisfies  
		$\EE_{\cC_\delta} \norm{\cC_\delta(\vx) -\vx}^2  \leq (1-\delta)\norm{\vx}^2, \forall \vx \in \RR^d$,  where $\delta\in(0,1]$.
	\end{assumption}
	\begin{remark}
		If $\delta=1$, there is no compression error, i.e., $\cC_\delta(\vx) =\vx$.  
		For instance, the Top-$k$ and Random-$k$ methods (see e.g., \cite{koloskova*2020decentralized,Beznosikov2020Biased}) satisfy Assumption \ref{Assumption:Biased}, where $\delta$ is given by $\delta=\frac{k}{p}$.
		
	\end{remark}
	
	\subsubsection{General compression operators}
	We now present a general assumption on the compression operators, which contains Assumptions \ref{Assumption:Unbiased} and \ref{Assumption:Biased} as special cases.
	\begin{assumption}\label{Assumption:General}
		The  compression operator $\cC \colon \RR^d \to \RR^d$   satisfies  
		\begin{align}\label{def:Ncompressor}
			\EE_{\cC} \norm{\cC(\vx) -\vx}^2 &\leq C\norm{\vx}^2,& ~\forall \vx \in \RR^d , 
		\end{align}
		and the $r$-scaling of $\cC$  satisfies 
		\begin{align}\label{def:contract}
			\EE_{\cC} \norm{\frac{\cC(\vx)}{r} -\vx}^2 &\leq (1-\delta)\norm{\vx}^2,& ~\forall \vx \in \RR^d , 
		\end{align}
		for some constants $\delta\in(0,1]$ and  $r>0$.
	\end{assumption}
	
	\begin{remark}
		On one hand, if $C<1$, Assumption \ref{Assumption:General} degenerates to Assumption \ref{Assumption:Biased} by setting $r=1$ and $\delta=1-C$. On the other hand, if $\cC$ is unbiased, i.e., $\EE \cC(\vx)=\vx$, then  Assumption \ref{Assumption:General} degenerates to Assumption  \ref{Assumption:Unbiased} by setting $r=C+1$ and $\delta=\frac{1}{C+1}$. In short, Assumption \ref{Assumption:General} gives a unified description of   unbiased and biased compression operators and thus Assumptions  \ref{Assumption:Unbiased} and \ref{Assumption:Biased} can be regarded as its special cases.
		
		However, there also exist compression operators where $\cC$ is biased and $C\geq 1$ in Assumption \ref{Assumption:General}, that is, they do not satisfy Assumptions \ref{Assumption:Unbiased} and  \ref{Assumption:Biased}. Examples include the norm-sign compressor $\cC(\vx)=\|\vx\|_{q} \text{sign}(\vx)$ and the composition of quantization and sparsification \cite{karimireddy2019Error,signSGD-ICML18,Basu2020qsparse}.
	\end{remark}
	
	\begin{remark}
		Although some compression operators (e.g.,   composition of quantization and sparsification) can be rescaled so that the new compression operator  satisfies the contractive condition in Assumption \ref{Assumption:Biased}, applying the rescaled operator may hurt the performance of the algorithm when compared with directly using the original compression operator $\cC$. Considering Assumption \ref{Assumption:General} provides us with more flexibility in choosing the most suitable compression method.
	\end{remark}
	
	\section{A Compressed Gradient Tracking Algorithm}\label{sec: Alg}
	In this section, we introduce the communication-efficient compressed gradient tracking algorithm (C-GT). We also give some interpretations as well as how C-GT connects to existing works.
	
	Denote $\vD=\text{diag}([\eta_1,\eta_2,\ldots,\eta_n])$, where $\eta_i$ is the step-size of agent $i$. Let \textbf{Compress} be the compression function, and the compression operators are associated with the function \textbf{Compress}.
	The proposed compressed gradient tracking  algorithm (C-GT) is presented in Algorithm \ref{Alg:CGTe}. 
	\begin{figure}[htp]
		\centering
		\begin{minipage}{.99\linewidth}
			\begin{algorithm}[H]
				\caption{A Compressed Gradient Tracking (C-GT) Algorithm}
				\label{Alg:CGTe}
				\textbf{Input:} stopping time $K$,  step-size $\{\eta_i\}$, consensus step-size $\gamma$, scaling parameters $\alpha_x, \alpha_y$, and initial values $\vX^0$, $\vH_{x}^0$, $\vH_{y}^0$, $\vY^0 = \nabla \vF(\vX^0)$\\ 
				\noindent\textbf{Output:} $\vX^K, \vY^K$
				\begin{algorithmic}[1]
					\State $\vH^0_{x,w} = \vW \vH_{x}^0$
					\State $\vH^0_{y,w} = \vW \vH_{y}^0$				
					\For{$k=0,1,2,\dots, K-1$}
					\State $\widehat{\vX}^k,\widehat{\vX}^k_w,\vH_{x}^{k+1},\vH_{x,w}^{k+1}=\Call{Comm}{\vX^k,\vH^k_x,\vH^k_{x,w}}$\footnote{In  Lines 4 and 5, $\alpha_z$ in the compression function is replaced by $\alpha_x$ for decision difference compression and $\alpha_y$ for gradient tracker difference compression, respectively.}
					\State $\widehat{\vY}^k,\widehat{\vY}^k_w,\vH_{y}^{k+1},\vH_{y,w}^{k+1}=\Call{Comm}{\vY^k,\vH^k_y,\vH^k_{y,w}}$
					\State $\vX^{k+1}=\vX^{k}-\gamma(\widehat{\vX}^{k}-\widehat{\vX}^k_w)-\vD \vY^{k}$
					\State $\vY^{k+1}=\vY^{k}-\gamma(\widehat{\vY}^{k}-\widehat{\vY}^k_w)+\nabla \vF(\vX^{k+1})-\nabla \vF(\vX^{k})$ 
					\EndFor
					\vspace{0.01in}
					\Procedure{Comm}{$\vZ, \vH, \vH_{w}$}
					\State $\vQ =$ \textbf{Compress}($\vZ - \vH$) \hfill $\vartriangleright$ Compression
					\State $\widehat \vZ = \vH + \vQ$
					\State $\widehat \vZ_w = \vH_{w} + \vW\vQ$ \hfill $\vartriangleright$ Communication
					\State $\vH \leftarrow (1-\alpha_z)\vH + \alpha_z \widehat \vZ $
					\State $\vH_{w} \leftarrow (1-\alpha_z)\vH_{w} + \alpha_z \widehat \vZ_w $
					\State \textbf{Return:} $\widehat \vZ, \widehat \vZ_w, \vH, \vH_{w}$
					\EndProcedure
				\end{algorithmic}
			\end{algorithm}
		\end{minipage}
	\end{figure}
	
	The compression and communication steps are included in the procedure $\Call{Comm}{\vZ, \vH, \vH_w}$.  The function \textbf{Compress} is the compression operator that independently compresses the variables for each agent per iteration. In  Line 10,  the difference between $\vZ$ and the auxiliary variable $\vH$ is compressed and then added back to $\vH$ in Line 11 for obtaining $\widehat{\vZ}$. Here $\vH$ acts as a reference point, and when it gradually approaches $\vZ$  such that the difference  vanishes to $0$,  the compression error on the difference will also decrease to $0$ under Assumption \ref{Assumption:General}. The low-bit compressed value $\vQ$ is transmitted in Line 12.
	
	To control the  compression error, particularly for a relatively large constant $C$ in Assumption \ref{Assumption:General}, we introduce a momentum update $\vH=(1-\alpha_z)\vH+\alpha_z \widehat \vZ$ motivated from the centralized distributed method DIANA \cite{mishchenko2019distributed} and the decentralized algorithm LEAD \cite{liu2020linear}. If $\alpha_z=1$, the update degenerates to that in the decentralized stochastic algorithm CHOCO-SGD \cite{ksj2019choco}.
	
	In Line 14, $\vH_w$ is used as a back-up copy for the neighboring information. 
	By introducing such an auxiliary variable, there is no need to store all the neighbors' reference points $\vH$ \cite{ksj2019choco}. 
	% 	This is a key step for ensuring the communication efficiency implementation since we do not need to store all the neighbors' reference points $\vH$. 
	Noticing that $\vH_w^{0}=\vW \vH^{0}$ from the initialization in Lines 1 and 2, we have $\widehat{\vZ}_w=\vW \widehat{\vZ}$ and then $\vH_w=\vW \vH$ by induction. It follows that $\widehat\vX_w^k=\vW\widehat\vX^k$ and $\widehat\vY_w^k=\vW\widehat\vY^k$ from Lines 4 and 5. Therefore, the decision variable update in Line 6 becomes
	\begin{align}\label{eq:EffGT_X}
		\nonumber\vX^{k+1}=&\vX^{k}-\gamma(\widehat{\vX}^{k}-\vW\widehat\vX^k)-\vD  \vY^{k}\\
		=&\vX^{k}-\gamma(\vI-\vW)\widehat{\vX}^{k}-\vD  \vY^{k},
	\end{align}
	and the gradient tracker update in Line 7 is given by
	\begin{multline}\label{eq:EffGT_Y}
		\vY^{k+1}=\vY^{k}-\gamma(\widehat{\vY}^{k}-\vW\widehat\vY^k)+\nabla \vF(\vX^{k+1})-\nabla \vF(\vX^{k})\\
		=\vY^{k}-\gamma(\vI-\vW)\widehat{\vY}^{k}+\nabla \vF(\vX^{k+1})-\nabla \vF(\vX^{k}).
	\end{multline}
	
	One key property of C-GT is that gradient tracking is efficient regardless of the compression errors, i.e., for $k\ge 0$,  
	\begin{align}\label{eq:avgY}
		\nonumber \vone^{\T}\vY^{k+1}=&\vone^{\T}(\vY^{k}-\gamma(\vI-\vW)\widehat{\vY}^{k}+\nabla \vF(\vX^{k+1})-\nabla \vF(\vX^{k}))\\
		\nonumber=&\vone^{\T}\vY^{k}+\vone^{\T}\nabla \vF(\vX^{k+1})-\vone^{\T}\nabla \vF(\vX^{k})\\
		=&\vone^{\T}\nabla \vF(\vX^{k+1}).
	\end{align}
	The second equality holds because $\vone^{\T}(\vI-\vW)=0$, and the last equality is obtained by induction under the initial condition $\vY^0 = \nabla \vF(\vX^0)$. Therefore, as long as $\vy_{i}^{k}$ reaches (approximate) consensus among all the agents, each $\vy_{i}^{k}$ is able to track the average gradient $\vone^{\T}\nabla \vF(\vX^{k})/n$.
	Moreover, by multiplying $\vone^{\T}$ and dividing $n$ on both sides of Line 8, we obtain
	\begin{align}\label{eq:avgX}
		\overline \vX^{k+1}=\overline \vX^{k}-\frac{1}{n}\vone^{\T}\vD \vY^{k}.
	\end{align}
	Hence the update of $\overline \vX^{k}$ does not involve any compression error from the current step.

	\begin{remark}
		If no communication compression is performed in the algorithm, i.e., $\widehat{\vX}^k=\vX^{k}$ and $\widehat{\vY}^k=\vY^{k}$,   
		then C-GT recovers the typical   distributed  gradient tracking algorithm in \cite{Nedic2017achieving,Qu2018Harnessing}  for $\eta_i=\eta$.  To see such a connection, note that C-GT reads
		\begin{align}\label{eq:GT_X}
			\nonumber\vX^{k+1}=&\vX^{k}-\gamma(\vI-\vW)\vX^{k}-\eta \vY^{k}\\
			=&[(1-\gamma)\vI+\gamma \vW] \vX^{k}-\eta \vY^{k},
		\end{align}
		and
		\begin{multline}\label{eq:GT_Y}
			\vY^{k+1}=\vY^{k}-\gamma(\vI-\vW)\vY^{k}+\nabla \vF(\vX^{k+1})-\nabla \vF(\vX^{k})\\
			=[(1-\gamma)\vI+\gamma \vW] \vY^{k}+\nabla \vF(\vX^{k+1})-\nabla \vF(\vX^{k}),
		\end{multline}
		where we substitute $\widehat{\vX}^k=\vX^{k}$ and $\widehat{\vY}^k=\vY^{k}$ in \eqref{eq:EffGT_X} and \eqref{eq:EffGT_Y}, respectively. By denoting $\widetilde{\vW}:=(1-\gamma)\vI+\gamma \vW$, C-GT takes the same form as the typical gradient tracking method.
	\end{remark}
	
	On the other hand, 	C-GT performs an implicit error compensation operation that mitigates the impact of  the compression error, as can be seen from the following argument. The decision variable is updated as
	\begin{align}\label{eq:error}
		\nonumber 	\vX^{k+1}=&\vX^{k}-\gamma(\vI-\vW)(\vX^{k}-\vE^{k})-\vD \vY^{k}\\
		=&[(1-\gamma)\vI+\gamma \vW] \vX^{k}  -\vD \vY^{k} + \gamma(\vI-\vW) \vE^{k},
	\end{align}
	where $\vE^{k}:=\vX^{k}-\widehat{\vX}^k$ measures the compression error for the decision variable. The additional  term $(\vI-\vW) \vE^k$ implies that each agent $i$ transmits its total compression error $-\sum_{j\in \mathcal{N}_{i}\cup \{i\} } w_{ji}\ve^k_i = -\ve^k_i$ to its neighboring agents and compensates this error locally by adding $\ve^k_i$, where $\ve^k_i\in\mathbb{R}^{1\times p}$ is the  $i$-th row of $\vE^{k}$.  
	Similarly, the compression errors for the gradient trackers are also mitigated.
	
	\section{Convergence Analysis}\label{sec: CA-CGT}
	In this section, we study the convergence properties of the
	proposed compressed gradient tracking algorithm for minimizing strongly convex and  smooth cost functions. 
	Our analysis relies on constructing a linear system of inequalities that is related to the
	optimization error $\Omega_{o}^{k}:=\EE\big[\|\oX^{k}-\vx^*\|^2\big]$, 
	consensus error $\Omega_{c}^{k}:=\EE\big[\|\vX^{k}-\vone\oX^{k}\|^2\big]$, 
	gradient tracking error $\Omega_{g}^{k}:=\EE\big[\|\vY^{k}-\vone\oY^{k}\|^2\big]$, 
	and compression errors $\Omega_{cx}^{k}:=\EE\big[\| \vX^{k}- \vH^{k}_x \|^2\big]$ and  $\Omega_{cy}^{k}:=\EE\big[\|\vY^{k}-\vH^{k}_y\|^2\big]$. 
	
	In order to derive the main results, we introduce some useful lemmas first.
	\begin{lemma}\label{lem1}
		Under Assumption \ref{Assumption: function},  for all $k\geq0$,  there holds
		\begin{equation}
			\|\nabla f((\oX^k)^{\T}) - (\oY^k)^{\T}\|\leq\frac{L}{\sqrt{n}}\|\vX^k-\vone\oX^k\|.
		\end{equation}
		In addition, if $\eta<2/(\mu+L),$ then we have
		\begin{equation}
			\|\vx-\eta\nabla f(\vx) - (\vx^*)^{\T}\|\leq(1-\eta\mu)\|\vx-(\vx^*)^{\T}\|,\forall \vx\in \RR^p.
		\end{equation}
	\end{lemma}
	\begin{lemma}\label{lem2}
		Suppose Assumption \ref{Assumption: network} holds. For any $\omega\in\RR^{n\times p}$, we have   
		$
		\|\vW\omega-\vone\bar{\omega}\|\leq\rho_w\|\omega-\vone\bar{\omega}\|,
		$ 
		where $\bar{\omega}=\frac{1}{n}\vone^\T\omega$ and $\rho_w<1$ is the spectral norm of the matrix $\vW-\frac{1}{n}\vone\vone^\T$. 
	\end{lemma}
	
	The proofs of Lemma \ref{lem1} and Lemma \ref{lem2} can be found in  Lemma 10 of \cite{Qu2018Harnessing}.
	\begin{remark}
		For  $\gamma \in (0,1]$, if we define $\widetilde{\vW}=(1-\gamma)\vI+\gamma \vW$, then 
		$
		\|\tW\omega-\vone\bar{\omega}\|\leq\tilde{\rho}\|\omega-\vone\bar{\omega}\|,
		$
		where $\tilde{\rho}  = 1-\gamma s,$ and $s = 1-\rho_w$.
	\end{remark}
	
	Denote $\bar{\eta}=\frac{1}{n}\sum_{i=1}^{n} \eta_i$, $\hat{\eta}=\max_{i} \eta_{i}$. We introduce below the key lemma for establishing the linear convergence of the C-GT algorithm under Assumptions \ref{Assumption: network}, \ref{Assumption: function} and \ref{Assumption:General}. 
	
	\begin{lemma}\label{Lem:CGT}
		Suppose Assumptions \ref{Assumption: network}, \ref{Assumption: function} and \ref{Assumption:General} hold and  $\hat{\eta}<\min\big\{\frac{2}{\mu+L},\frac{1}{3\mu}\big\}$. 
		Then we have the following linear system of inequalities:
		$$
		\vw^{k+1}\leq \vA\vw^k,
		$$
		where 
		$
		\vw^k :=	\big[\Omega_{o}^{k}, 
		\Omega_{c}^{k}, 
		\Omega_{g}^{k}, 
		\Omega_{cx}^{k}, 
		\Omega_{cy}^{k}\big]^{\T}.$ 
		The inequality is to be taken component-wisely, where $\vA\in\mathbb{R}^{5\times 5}$ is non-negative.\footnote{The elements of the transition matrix $\vA$ correspond to the parameters of the inequalities in the proof.}
	\end{lemma}
	\begin{IEEEproof}
		See Appendix \ref{Pf:LemCGT}.
	\end{IEEEproof}
	
	Based on Lemma \ref{Lem:CGT}, we present the preliminary convergence result for the C-GT algorithm below. 
	\begin{lemma}\label{Thm:CGT}
		Suppose Assumptions \ref{Assumption: network}, \ref{Assumption: function} and \ref{Assumption:General} hold,   the scaling parameters $\alpha_x, \alpha_y \in (0, \frac{1}{r}]$, $\bar{\eta}\geq M \hat{\eta}$ for some $M>0$, 
		and 
		\begin{align}
			&n\epsilon_1 \geq  \frac{12\kappa^2}{M^2}\epsilon_3,~
			\epsilon_2\geq \frac{12C\lambda}{s^2}\epsilon_4,~
			\epsilon_3\geq M\epsilon_2,~ \nonumber \\
			&\epsilon_3\geq \frac{12\lambda\left(3\epsilon_2+C(3\epsilon_4+\epsilon_5)\right)}{s^2},~
			\epsilon_4>0,~
			\epsilon_5>0,~ \label{epsilon1_5}\\ 
			&\gamma\leq \min\Big\{1,
			\frac{1-c_x}{m_x}\epsilon_4,
			\frac{1-c_y}{m_y}\epsilon_5\Big\}, \label{gamma_step}\\ 
			&\hat{\eta}\le \min\Big\{ 
			\sqrt{\frac{\epsilon_2}{12}},
			\sqrt{\frac{\epsilon_3}{36}}
			\Big\}\frac{s\gamma}{\sqrt{2n\epsilon_1+2\epsilon_2+\epsilon_3}L},\label{eta_step}
		\end{align}
		where $s=1-\rho_w$, $\lambda:=\norm{\vI-\vW}^2$, 
		$m_x:=t_x(2n\epsilon_1+2\epsilon_2+\epsilon_3)+t_x\lambda(\epsilon_2+C\epsilon_4)
		+\frac{M\epsilon_4}{2\kappa}
		$, 
		$
		m_y:=3t_y(2n\epsilon_1+2\epsilon_2+\epsilon_3)
		+t_y\lambda(3\epsilon_2+\epsilon_3
		+3C\epsilon_4+C\epsilon_5)
		+\frac{M\epsilon_5}{2\kappa}
		$,  $c_x=\tau_x(1-\alpha_x r \delta)<1$, $c_y=\tau_y(1-\alpha_y r \delta)<1$, $t_x=\frac{3\tau_x}{\tau_x-1}$, 
		$t_y=\frac{3\tau_y}{\tau_y-1}$, constants $\tau_x,\tau_y>1$,
		% 			are constants that are grater than one,
		and $\epsilon_1$-$\epsilon_5$ are some positive constants. 
		Then the spectral radius of $\vA$ satisfies $\rho\big(\vA\big)\leq 1-\frac{1}{2}M\hat{\eta}\mu$, and the optimization error $\Omega_{o}^{k}$ and the consensus error $\Omega_{c}^{k}$ both converge to 0 at the linear rate $\cO((1-\frac{1}{2}M\hat{\eta}\mu)^k)$.
	\end{lemma}
	\begin{IEEEproof}
		See Appendix \ref{Pf:ThmCGT}.
	\end{IEEEproof}

	\subsection{Main Results}
	By taking some concrete values for the constants in Lemma \ref{Thm:CGT}, we derive the main convergence result for the C-GT algorithm under Assumptions \ref{Assumption: network}, \ref{Assumption: function} and \ref{Assumption:General} 
	in the following theorem, which demonstrates the linearly convergent property of C-GT for the general compression operators. 
	
	\begin{theorem}\label{Thm1:C-GT}
		Suppose Assumptions \ref{Assumption: network}, \ref{Assumption: function}, \ref{Assumption:General} hold,  $\bar{\eta}\geq M \hat{\eta}$ for some $M>0$,   $\alpha_x,\alpha_y=1/r$, the consensus step-size $\gamma$ satisfies 
		\begin{align}
			\gamma\leq \min\Big\{1,
			\frac{\delta}{m(2-\delta)}\Big\},
		\end{align} 
		and the maximum step-size $\hat{\eta}$
		satisfies 
		\begin{align}
			\hat{\eta}\leq  \frac{s^2}{\sqrt{12\left[\left(\frac{24\kappa^2}{M^2}+1\right)(4s^2+36\lambda)+2s^2\right]}}
			\frac{\gamma}{L},
		\end{align}
		where $\kappa=L/\mu$ and
		$m=\frac{6}{\delta}\Big[(\frac{72\kappa^2}{M^2}+\lambda+3)(\frac{432C\lambda^2}{s^4}+\frac{48C\lambda}{s^2}) +(3\lambda+6)\frac{12C\lambda}{s^2}+4C\lambda\Big]+\frac{M}{2\kappa}$ with $\lambda=\norm{\vI-\vW}^2$. 
		Then, the optimization error $\Omega_{o}^{k}$ and the consensus error $\Omega_{c}^{k}$ both converge to 0 at the linear rate $\cO((1-\frac{1}{2}M\hat{\eta}\mu)^k)$.
	\end{theorem}
	\begin{IEEEproof}
		See Appendix \ref{Pf:Thm1CGT}.
	\end{IEEEproof}

	\begin{remark}
		If the  compression error is sufficiently small, i.e., $C\rightarrow0$ and $\delta\rightarrow1$, we have  $m\approx\frac{M}{2\kappa}<1$.
		Then, we obtain  $\gamma\leq 1$ and   
		\begin{align*}
			\hat{\eta}\leq & \frac{s^2}{\sqrt{12\left[\left(\frac{24\kappa^2}{M^2}+1\right)(4s^2+36\lambda)+2s^2\right]}}
			\frac{1}{L} {\sp\sim\mathcal{O}\left(\frac{(1-\rho_w)^2}{\kappa L}\right)}.
		\end{align*}
		The convergence rate of C-GT is then comparable to those of the typical gradient tracking methods; see, e.g., \cite{Qu2018Harnessing}.
	\end{remark}

	\begin{remark}
		In practice, the restrictions on $\alpha_x$ and $\alpha_y$ can be relaxed to $\alpha_x,\alpha_y\in(0,\frac{1}{r}]$ as in Lemma \ref{Thm:CGT}. 
		The condition $\bar{\eta}\geq M \hat{\eta}$ is always satisfied for some fixed $M$, e.g., $M=\frac{1}{n}$. If in addition that all $\eta_i$ are equal, then we can take $M=1$. 
	\end{remark}

	\begin{remark}
		Comparing the performance of C-GT with the existing linearly convergent algorithm LEAD \cite{li2019decentralized}, C-GT enjoys more flexibility in the mixing matrix, compression methods and the stepsize policy, while LEAD achieves faster convergence in theory under more restricted conditions.
	\end{remark}
	
	\section{Numerical Examples}\label{sec: simulation}
	In this part, we provide some numerical examples to confirm our theoretical results.  
	Consider the ridge regression problem:
	\begin{align}\label{Ridge Regression}
		\min_{x\in \mathbb{R}^{p}}f(x)=\frac{1}{n}\sum_{i=1}^nf_i(x)\left(=\left(u_i^{\T} x-v_i\right)^2+\rho\|x\|^2\right),
	\end{align}
	where $\rho>0$ is a penalty parameter. The pair $(u_i,v_i)$ is a sample that belongs to the $i$-th agent, where $u_i\in\mathbb{R}^p$ represents the features and $v_i\in\mathbb{R}$ represents the observations or  outputs.  
	In the simulations, pairs $(u_i,v_i)$ are pre-generated: input $u_i\in[-1,1]^p$ is uniformly distributed, and the output $v_i$ satisfies $v_i=u_i^{\T} \tilde{x}_i+\varepsilon_i$, where $\varepsilon_i$ are independent Gaussian noises with mean $0$ and variance $25$, and $\tilde{x}_i$ are predefined parameters evenly located in $[0,1]^p$. 
	Then, the $i$-th agent can calculate the gradient of its local objective function $f_i(x)$ with $g_i(x,u_i,v_i)=2(u_i^{\T}x -v_i)u_i+2\rho x$. The unique optimal solution of the problem is $x^*=(\sum_{i=1}^n u_i u_i^{\T}+n\rho\mathbf{I})^{-1}\sum_{i=1}^n u_i v_i$. 
	
	In our experimental settings, we consider penalty parameter $\rho=0.1$. The number of nodes is $n=100$, and the dimension of variables is $p=500$. Meanwhile, $\vx_i^0$ is randomly generated in $[0,1]^p$ and other initial values   satisfy $\vH_{x}^0=\vzero$, $\vH_{y}^0=\vzero$, and $\vY^0 = \nabla \vF(\vX^0)$.

	We   compare C-GT with CHOCO-SGD \cite{ksj2019choco}, LB \cite{kovalev2021linearly}, LEAD \cite{liu2020linear}  and the uncompressed linearly convergent methods, NIDS\cite{li2019decentralized} and GT\cite{Qu2018Harnessing}, for decentralized optimization over a randomly generated  undirected graph. In order to guarantee the fairness, all algorithms use their equivalent matrix forms. 
	The considered compression methods are 2-bit $\infty$-norm quantization (Q), Top-10 sparsification, composition of quantization and sparsification (Q-T) and its rescaled version (Q-T-R). Note that Q-T only satisfies Assumption \ref{Assumption:General} and does not satisfy Assumptions \ref{Assumption:Unbiased} and \ref{Assumption:Biased}. The communication bits of these compression methods are given in \cite{ksj2019choco,koloskova*2020decentralized,liu2020linear,Basu2020qsparse}.\footnote{ Note that C-GT requires two times the communication bits of the other algorithms per iteration since it compressed both the decision variable and the gradient tracker.} The parameter settings of the algorithms are given in Table \ref{table:parametersetUnd}, which are hand-tuned to achieve the best performance for each algorithm.

\begin{figure}[htp]
		\begin{center}
			\vspace{-1em}
			\includegraphics[width=7cm]{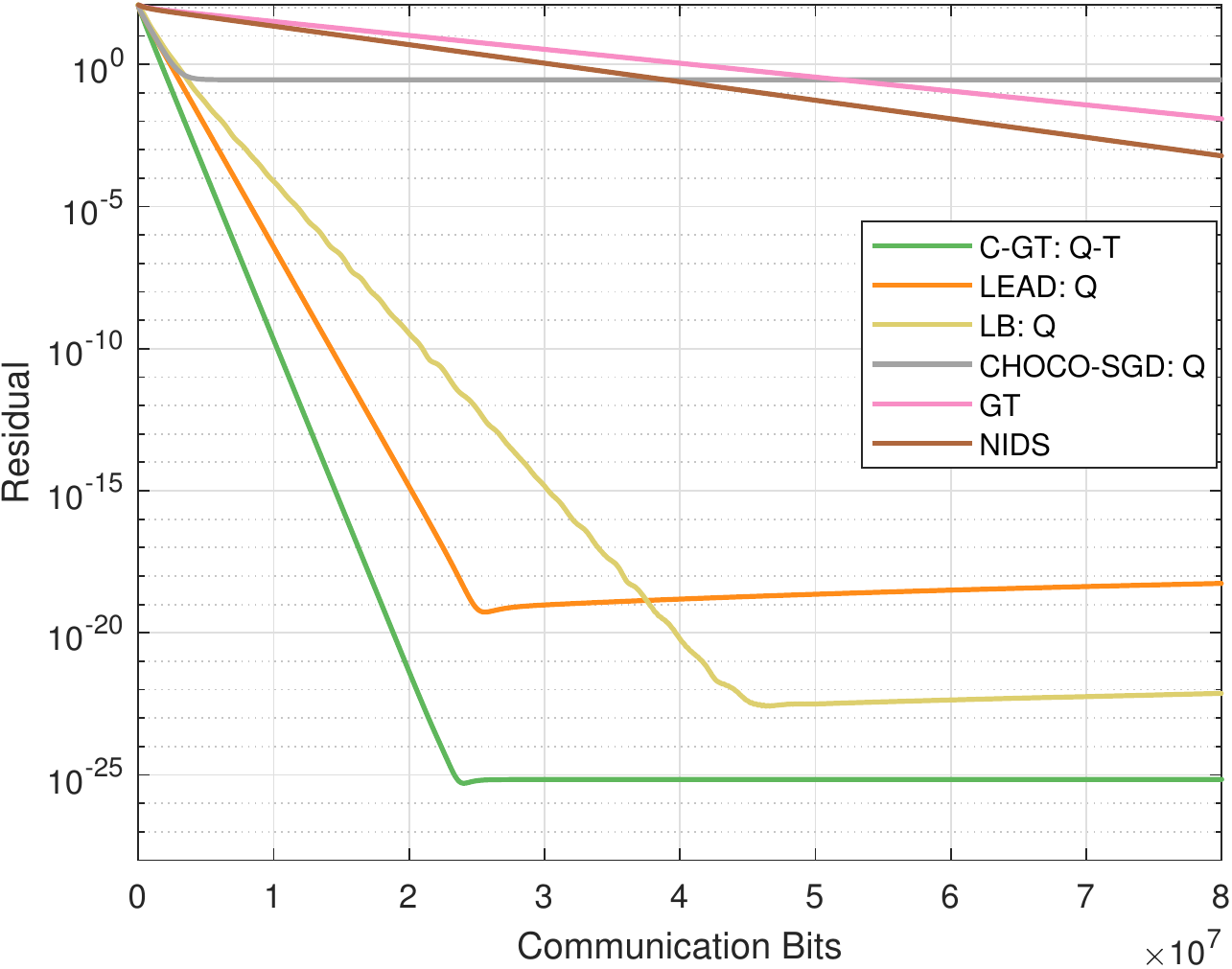}\vspace{-1em}
			\caption{ Residuals $\EE\big[\|\oX^{k}-\vx^*\|^2\big]$  against the communication bits for C-GT, LEAD, LB,  CHOCO-SGD and the uncompressed methods GT and NIDS.}\label{Fig1}
		\end{center}
	\end{figure}

	In Fig. \ref{Fig1}, we compare the communication efficiency of C-GT, LEAD, LB and CHOCO-SGD with the uncompressed methods GT and NIDS. For C-GT, LEAD, LB and CHOCO-SGD, we apply the compressors that work the best for them, respectively. Apparently, C-GT and LEAD outperform the other methods, while C-GT achieves the best communication efficiency.
	
\begin{figure}[htp]
		\begin{center}
			\vspace{-1em}
			\includegraphics[width=7cm]{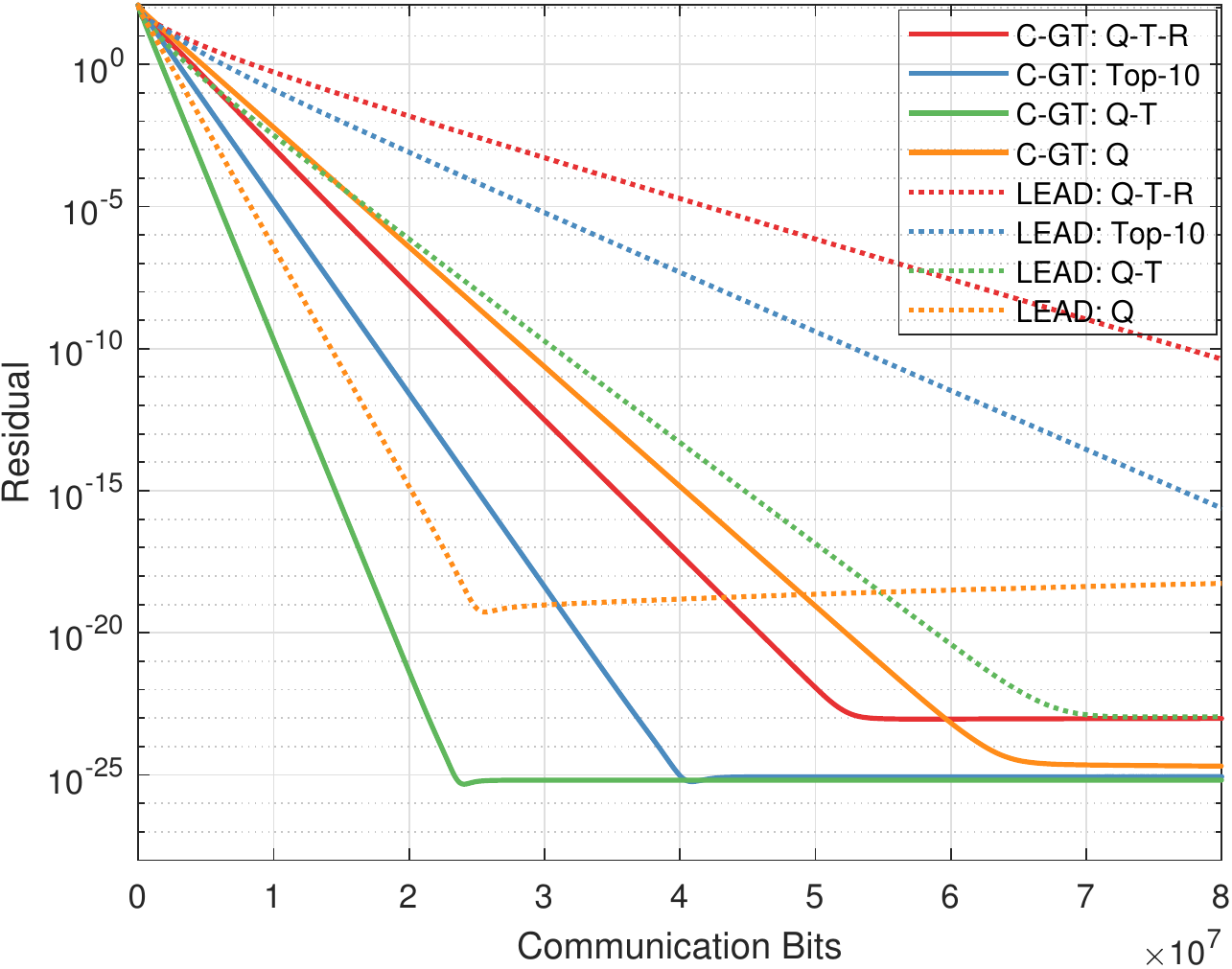}\vspace{-1em}
			\caption{ Residuals $\EE\big[\|\oX^{k}-\vx^*\|^2\big]$ against the communication bits for C-GT and LEAD under different compression methods.}\label{Fig2}
		\end{center}
	\end{figure}

	In Fig. \ref{Fig2}, we further present a detailed comparison between C-GT and LEAD under different types of compressors.
	Note that LEAD works the best under the unbiased compressor Q, while C-GT is more efficient under the biased compression operator Top-10 and the composition of quantization and sparsification Q-T. In particular, the performance of C-GT under Q-T is the most favorite among all the combinations. 
	
	It can also be seen from Fig. \ref{Fig2} that using the rescaled compressors
	leads to slower convergence, which suggests
	that rescaling the compression operators to satisfy the typical
	contractive requirement (i.e., Assumption \ref{Assumption:Biased}) may harm the
	algorithmic performance. Therefore, we can conclude that considering Assumption \ref{Assumption:General} provides
	users with more freedom in choosing the best compression
	method. These experimental findings demonstrate the effectiveness of C-GT.

	\begin{table}[htp]
		\begin{center}
			\begin{tabular}{cccccc}
				\hline
				Algorithm       &Compressor   &$\alpha_x$&$\alpha_y$  &$\gamma$       &$\eta$\\
				\hline
				LEAD            &Top-10        &1         &1   &$7\times 10^{-5}$    &0.005\\
				C-GT            &Top-10        &1         &1           &0.05         &0.005\\
				LEAD            &Q             &1         &1           &$0.09$       &0.005\\
				C-GT            &Q             &1         &1           &$0.09$       &0.005\\
				CHOCO-SGD       &Q           &1         &1           &$0.7$  &0.005\\	
				LB       &Q           &1         &1           &$0.1$  &0.005\\	
				LEAD            &Q-T           &1         &1   &$7\times 10^{-5}$    &0.005\\
				C-GT            &Q-T           &1         &1           &$0.06$       &0.005\\
				LEAD            &Q-T-R         &1         &1   &$3\times 10^{-5}$    &0.005\\
				C-GT            &Q-T-R         &1         &1           &0.009  &0.005\\		
				GT              &/             &/         &/           &/  &0.009\\	
				NIDS            &/             &/         &/           &/  &0.006\\			
				\hline
			\end{tabular}
		\end{center}
		\caption{Parameter setting  for different algorithms and compression methods.}	
		\label{table:parametersetUnd}
	\end{table}
	
	
	\section{Conclusions}\label{sec: conclusion}
	In this paper, we consider the problem of decentralized optimization with communication compression over a multi-agent network. 
	Specifically, we   propose  a compressed gradient tracking algorithm, termed C-GT, and show the algorithm converges linearly for strongly convex and smooth objective functions. 
	C-GT not only inherits the advantages of gradient tracking-based methods, but also works with a wide class of compression operators.  
	Simulation examples demonstrate the  effectiveness and flexibility of C-GT for  undirected networks. 
	Future work will consider equipping C-GT with accelerated techniques such as Nesterov's acceleration and momentum methods. Non-convex objective functions are also of future concern.
	%%%%%%%%%%%%%%%%%%%%%%%%%%%%%%%%%%%%%%%%%%%%%%%%%%%%%%%%%%%%%%%%%%%%%%%%%%%%%%%%%%%%%%%%%%%%%%%%%%	
	\section*{ACKNOWLEDGMENT}
	We would like to thank Bingyu Wang from The Chinese University of Hong Kong, Shenzhen for providing helpful feedback.
	%%%%%%%%%%%%%%%%%%%%%%%%%%%%%%%%%%%%%%%%%%%%%%%%%%%%%%%%%%%%%%%%%%%%%%%%%%%%%%%%%%%%%%%%%%%%%%%%%
	%%%%%%%%%%%%%%%%%%%%%%%%%%%%%%%%%%%%%%%%%%%%%%%%%%%%%%%%%%%%%%%%%%%%%%%%%%%%%%%%%%%%%%%%%%%%%%%%%%	
	\bibliographystyle{IEEEtran}
	\bibliography{lywbib_compression}	
	%%%%%%%%%%%%%%%%%%%%%%%%%%%%%%%%%%%%%%%%%%%%%%%%%%%%%%%%%%%%%%%%%%%%%%%%%%%%%%%%%%%%%%%%%%%%%%%%%%
	%\addtolength{\textheight}{-12cm}   
	% This command serves to balance the column lengths
	% on the last page of the document manually. It shortens
	% the textheight of the last page by a suitable amount.
	% This command does not take effect until the next page
	% so it should come on the page before the last. Make
	% sure that you do not shorten the textheight too much.
	
	%\appendices
	\appendix
	
	\section{Proofs for the C-GT Algorithm}\label{AppenddixB}
	\subsection{Supplementary Lemmas}
	The following vector and matrix inequalities are often invoked.
	\begin{lemma}\label{lem:UV}
		For  $\vU,\vV\in \RR^{n\times p}$ and any constant $\tau>0$,   we have the following inequality:
		\begin{align}
			\|\vU+\vV\|^2\leq (1+\tau)\|\vU\|^2 + (1+\frac{1}{\tau})\|\vV\|^2.
		\end{align}
		In particular, taking $\tau=\tau'-1, \tau'>1$, we have
		\begin{align}\label{ineq:UV}
			\|\vU+\vV\|^2\leq \tau'\|\vU\|^2 + \frac{\tau'}{\tau'-1}\|\vV\|^2.
		\end{align}
		
		In addition, for any $\vU_1,\vU_2,\vU_3 \in \RR^{n\times p}$, we have 
		$
		\|\vU_1+\vU_2+\vU_3\|^2\leq \tau'\|\vU_1\|^2 + \frac{2\tau'}{\tau'-1}\left[\|\vU_2\|^2+\|\vU_3\|^2\right]
		$ 	
		and 
		$
		\|\vU_1+\vU_2+\vU_3\|^2\leq 3\|\vU_1\|^2 + 3\|\vU_2\|^2+3\|\vU_3\|^2.
		$
	\end{lemma}
	
	\begin{lemma}\label{implma}
		(Corollary 8.1.29 in \cite{Horn2012Matrix}) Let $\vM\in\RR^{l\times l}$ and $\vv\in\RR^l$ be a nonnegative matrix  and an element-wise  positive vector, respectively. If $\vM\vv\le \theta \vv$, then $\rho(\vM)\le \theta.$
	\end{lemma}

	\subsection{Proof of Lemma \ref{Lem:CGT}}\label{Pf:LemCGT}
	Let $\cF^k$ be the $\sigma$-algebra generated by $\{\vX^0,\vY^0,\vX^1,\vY^1,\cdots,\vX^{k},\vY^{k}\}$, and define $\EE[ \cdot |\cF^k]$ as the conditional expectation with respect to the compression operator given $\cF^k$.
	
	Before deriving the linear system of inequalities, we bound $\EE\left[\|\vX^k-\widehat{\vX}^k\|^2\middle |\cF^k\right]$ and $\EE\left[\|\vY^k-\widehat{\vY}^k\|^2\middle |\cF^k\right]$, respectively. 
	From the $\Call{COMM}{}$ procedure in Algorithm \ref{Alg:CGTe} for the decision variable, we know $\norm{\vX^k-\widehat{\vX}^k}^2=\norm{\vX^k-\vH^k_x-\cC(\vX^k-\vH^k_x)}^2$. Then, we obtain
	\begin{align}\label{ineqN:Xk_Xkhat} 
		\EE\left[\norm{\vX^k-\widehat{\vX}^k}^2\middle |\cF^k\right]
		\leq C\norm{\vX^k-\vH^k_x}^2.
	\end{align}
	Similarly, we have
	\begin{align}\label{ineqN:Yk_Ykhat}
		\EE\left[\norm{\vY^k-\widehat{\vY}^k}^2\middle |\cF^k\right]\leq C\norm{\vY^k-\vH^k_y}^2.
	\end{align}
	
	\subsubsection{Optimality error}
	According to  \eqref{eq:avgX},  Lemmas \ref{lem1} and \ref{lem:UV}, we obtain
	\begin{align*}
		& \norm{\oX^{k+1} - \vx^*}^2 \\
		=& \Big\|\oX^k-\bar{\eta}\nabla f((\oX^k)^\T)^\T-\vx^*\\
		&\quad-\bar{\eta}(\oY^k-\nabla f((\oX^k)^\T)^\T)-\frac{1}{n}\vone^{\T}(\vD-\bar{\eta}\vI) (\vY^k-\vone\oY^k)\Big\|^2\\
		\leq & (1+\tau_1)(1-\bar{\eta}\mu)^2 \norm{\oX^k-\vx^*}^2\\
		&+2(1+\frac{1}{\tau_1})\Big(\frac{\bar{\eta}^2 L^2}{n}\norm{\vX^k-\vone \oX^k}^2\\
		&\qquad+\frac{1}{n}\norm{\vD-\bar{\eta}\vI}^2\norm{\vY^k-\vone\oY^k}^2\Big).
	\end{align*}
	Taking $\tau_1=\frac{3}{8}\bar{\eta}\mu$ and noticing that $\bar{\eta}\leq\hat{\eta}\leq\frac{1}{3\mu}$, we have
	\begin{align*}
		&\norm{\oX^{k+1} - \vx^*}^2\\
		\leq&(1+\frac{3}{8}\bar{\eta}\mu)(1-\bar{\eta}\mu)^2\norm{\oX^k - \vx^*}^2 \\
		&+ 2(1+\frac{8}{3\bar{\eta}\mu})\Big(\frac{\bar{\eta}^2 L^2}{n}\norm{\vX^k-\vone \oX^k}^2+\frac{\hat{\eta}^2}{n}\norm{\vY^k-\vone\oY^k}^2\Big)\\
		\leq&(1-\frac{3}{2}\bar{\eta}\mu) \norm{\oX^k-\vx^*}^2 \\
		&+ \frac{6\bar{\eta} L^2}{\mu n}\norm{\vX^k-\vone \oX^k}^2+\frac{6\hat{\eta}^2}{\bar{\eta}n\mu}\norm{\vY^k-\vone\oY^k}^2\\
		\leq&(1-\frac{3}{2}M\hat{\eta}\mu) \norm{\oX^k-\vx^*}^2 \\
		&+ \frac{6\hat{\eta} L^2}{\mu n}\norm{\vX^k-\vone \oX^k}^2+\frac{6\hat{\eta}}{\mu nM}\norm{\vY^k-\vone\oY^k}^2.
	\end{align*}
	In the first inequality, we use $\norm{\vD-\bar{\eta}\vI}\leq \hat{\eta}$. 
	The last inequality holds since $\bar{\eta}\leq\hat{\eta}$ and $\bar{\eta}\geq M\hat{\eta}$ for some $M>0$.			
	\subsubsection{Consensus error}
	Recalling   relations \eqref{eq:EffGT_X} and  \eqref{eq:avgX}, we know
	\begin{align*}
		&\EE \left[\norm{\vX^{k+1} - \vone \oX^{k+1}}^2\middle|\cF^k\right] \\
		=&\EE \Big[\Big \|\vX^k - \gamma(\vI-\vW)\widehat{\vX}^k -\vD \vY^k\\
		&-\vone (\oX^k - \frac{1}{n}\vone^\T \vD \vY^k) \Big\|^2 \bigg| \cF^k\Big]\\
		=&\EE \bigg[\Big\|\gamma(\vI-\vW)(\vX^k-\widehat{\vX}^k) +\vX^k -\gamma(\vI-\vW)\vX^k -\vone \oX^k\\
		&-(\vI-\frac{1}{n}\vone\vone^\T) \vD \vY^k\Big\|^2 \bigg| \cF^k\bigg]\\
		=&\EE \bigg[\Big\|(\tW\vX^k-\vone\oX^k)-(\vI-\frac{1}{n}\vone\vone^\T) \vD \vY^k \\
		&+\gamma(\vI-\vW)(\vX^k-\widehat{\vX}^k)\Big\|^2 \bigg| \cF^k \bigg].
	\end{align*}
	Based on Lemmas \ref{lem2} and \ref{lem:UV}, we have 
	\begin{align*}
		&\EE \left[\norm{\vX^{k+1} - \vone \oX^{k+1}}^2\middle|\cF^k\right] \\
		\leq&(1+\tau_2)\tilde{\rho}^2\norm{\vX^k-\vone\oX^k}^2 
		+2(1+\frac{1}{\tau_2})\bigg(\norm{(\vI-\frac{1}{n}\vone\vone^\T)\vD\vY^k}^2\\
		&+\EE\bigg[\norm{\gamma(\vI-\vW)(\vX^k-\widehat{\vX}^k)}^2 \bigg| \cF^k\bigg]\bigg),
	\end{align*}
	where $\tilde{\rho}=1-s\gamma$ and $s=1-\rho_w$. Taking $\tau_2 = \frac{1-\tilde{\rho}^2}{2 \tilde{\rho}^2}$, we obtain
	\begin{align*}
		&\EE \left[\norm{\vX^{k+1} - \vone \oX^{k+1}}^2\middle|\cF^k\right] \\
		\leq &\frac{1+\tilde{\rho}^2}{2}\norm{\vX^k-\vone\oX^k}^2 +2\frac{1+\tilde{\rho}^2}{1-\tilde{\rho}^2}\Big(\norm{(\vI-\frac{1}{n}\vone\vone^\T) \vD\vY^k}^2\\
		&+\gamma^2\norm{\vI-\vW}^2\EE \left[\norm{\vX^k-\widehat{\vX}^k}^2 | \cF^k\right]\Big).
	\end{align*}
	Noticing that $\norm{\vI-\frac{1}{n}\vone\vone^\T}=1$ and  
	$
	\frac{1+\tilde{\rho}^2}{1-\tilde{\rho}^2}
	=\frac{1}{1-\tilde{\rho}}\frac{1+\tilde{\rho}^2}{1+\tilde{\rho}}
	\leq \frac{1}{1-\tilde{\rho}}
	=\frac{1}{s\gamma}
	$, 
	we derive
	\begin{align}\label{ineq:AvgConsensus}
		\nonumber 	&\EE \left[\norm{\vX^{k+1} - \vone \oX^{k+1}}^2\middle|\cF^k\right] \\
		\nonumber 	\leq &\frac{1+\tilde{\rho}^2}{2}\norm{\vX^k-\vone\oX^k}^2 +\frac{2\hat{\eta}^2}{s\gamma}\norm{\vY^k}^2\\
		&+\frac{2\gamma}{s}\norm{\vI-\vW}^2 \EE \left[\norm{\vX^k-\widehat{\vX}^k}^2 \middle| \cF^k\right]
	\end{align}
	In addition, we have 
	\begin{align}%\label{Yk}
		\nonumber \norm{\vY^k}^2=&\norm{\vY^k-\vone\oY^k+\vone\oY^k}^2\\
		\label{eq:Yk}=&\norm{\vY^k-\vone\oY^k}^2+n\norm{\oY^k}^2,
	\end{align}
	and 
	\begin{align}%\label{Ykbar}
		\nonumber \norm{\oY^k}^2=&\norm{\frac{1}{n}\sum_{i=1}^{n}\nabla f_i((\vx_i^k)^\T)-\nabla f((\oX^k)^\T)+\nabla f((\oX^k)^\T)^\T}^2\\
		\nonumber \leq&2\norm{\frac{1}{n}\sum_{i=1}^{n}[\nabla f_i((\vx_i^k)^\T)-\nabla f_i((\oX^k)^\T)]}^2\\
		\nonumber &+2\norm{\nabla f((\oX^k)^\T)^\T-\nabla f((\vx^*)^\T)^\T}^2\\
		\leq& \frac{2L^2}{n}\norm{\vX^k-\vone\oX^k}^2+2L^2\norm{\oX^k-\vx^*}^2,
		\label{ineq:Ykbar}
	\end{align}
	where the last inequality stems from the convexity of 2-norm and Assumption \ref{Assumption: function}.
	Plugging  \eqref{ineqN:Xk_Xkhat}, \eqref{eq:Yk} and \eqref{ineq:Ykbar} into \eqref{ineq:AvgConsensus}, we get
	\begin{align*}
		&\EE \left[\norm{\vX^{k+1} - \vone \oX^{k+1}}^2\middle|\cF^k\right] \\
		\leq&\Big(\frac{1+\tilde{\rho}^2}{2}+\frac{4L^2\hat{\eta}^2 }{s\gamma}\Big)\norm{\vX^k-\vone\oX^k}^2 +\frac{2\hat{\eta}^2}{s\gamma}\norm{\vY^k-\vone \oY^k}^2\\
		&+ \frac{2C\lambda\gamma}{s}\norm{\vX^k-\vH_x^k}^2+ \frac{4nL^2\hat{\eta}^2  }{s\gamma}\norm{\oX^k-\vx^{*}}^2 ,
	\end{align*}
	where $\lambda=\norm{\vI-\vW}^2$.
	\subsubsection{Gradient tracker error}
	For notational convenience, we use $\nabla^k:=\nabla \vF(\vX^{k})$ instead. 	Firstly, from relation \eqref{eq:EffGT_Y} and \eqref{eq:avgY}, the gradient tracker error is bounded by  
	\begin{align*}
		&\EE \left[\norm{\vY^{k+1} - \vone \oY^{k+1}}^2\middle|\cF^k\right] \\
		=& \EE \Big[\Big\|\vY^{k}-\gamma(\vI-\vW)\widehat{\vY}^k+\nabla^{k+1}-\nabla^k\\
		&\qquad -\vone (\oY^k+\frac{1}{n}\vone^\T\nabla^{k+1}-\frac{1}{n}\vone^\T\nabla^{k})\Big\|^2 \bigg|\cF^k\Big] \\
		=& \EE \bigg[\Big\|(\tW\vY^k-\vone\oY^k)+\gamma(\vI-\vW)(\vY^{k}-\widehat{\vY}^k)\\
		&+(\vI-\frac{1}{n}\vone\vone^\T)(\nabla^{k+1}-\nabla^k)\Big\|^2\bigg|\cF^k\bigg] \\
		\leq& (1+\tau_3)\tilde{\rho}^2\norm{\vY^k-\vone\oY^k}^2+2(1+\frac{1}{\tau_3})\EE \bigg[\norm{\nabla^{k+1}-\nabla^k}^2\\
		&+\norm{\gamma(\vI-\vW)(\vY^k-\widehat{\vY}^k)}^2 \bigg| \cF^k\bigg]
	\end{align*}
	Similar to the derivation of \eqref{ineq:AvgConsensus}, we obtain 
	\begin{align}
		\nonumber &\EE \left[\norm{\vY^{k+1} - \vone \oY^{k+1}}^2\middle|\cF^k\right] \\
		\nonumber \leq&\frac{1+\tilde{\rho}^2}{2}\norm{\vY^k-\vone\oY^k}^2 +\frac{2}{s\gamma}\EE\left[\norm{\nabla^{k+1}-\nabla^k}^2\right]\\
		\label{YkplusCGT}
		&+\frac{2\gamma}{s}\norm{\vI-\vW}^2\EE\left[\norm{\vY^k-\widehat{\vY}^k}^2 \bigg| \cF^k\right].
	\end{align}
	% 			Next, we derive the bound
	Given that
	\begin{align*}%\label{GradBound}
		\EE\left[\norm{\nabla^{k+1}-\nabla^k}^2\middle| \cF^k\right]\leq L^2\EE\left[\norm{\vX^{k+1}-\vX^k}^2\middle| \cF^k\right],
	\end{align*}
	we need to bound 
	\begin{align}
		\nonumber &\EE\left[\norm{\vX^{k+1}-\vX^k}^2\middle| \cF^k\right]\\
		\nonumber =&\EE\bigg[\Big\|\gamma(\vI-\vW)(\vX^k-\widehat{\vX}^k)\\
		\nonumber &\qquad-\gamma(\vI-\vW)(\vX^k-\vone\oX^k)-\vD\vY^k\Big\|^2\bigg| \cF^k\bigg]\\
		\nonumber\leq&3\gamma^2\norm{\vI-\vW}^2\EE\left[\norm{\vX^k-\widehat{\vX}^k}^2\middle| \cF^k\right]\\
		\label{ineq:XkplusminusXk}
		&+3\gamma^2\norm{\vI-\vW}^2\norm{\vX^k-\vone\oX^k}^2+3\hat{\eta}^2\norm{\vY^k}^2.
	\end{align}
	Plugging \eqref{ineqN:Xk_Xkhat}, \eqref{eq:Yk} and \eqref{ineq:Ykbar} into  \eqref{ineq:XkplusminusXk}, we obtain
	\begin{multline}
		\EE\left[\norm{\vX^{k+1}-\vX^k}^2\middle| \cF^k\right]
		\leq 6nL^2\hat{\eta}^2\norm{\oX^k-\vx^*}^2\\
		+(6L^2\hat{\eta}^2+3\gamma^2\norm{\vI-\vW}^2)\norm{\vX^k-\vone\oX^k}^2
		+3\hat{\eta}^2\norm{\vY^k-\vone\oY^k}^2\\
		\label{ineq:XkplusminusXk2}+3C\gamma^2\norm{\vI-\vW}^2\norm{\vX^k-\vH^k_x}^2.
	\end{multline}
	Therefore, we get
	\begin{multline}\label{GradplusminusGrad}
		\EE\left[\norm{\nabla^{k+1}-\nabla^k}^2\middle| \cF^k\right]
		\leq  6nL^4\hat{\eta}^2\norm{\oX^k-\vx^*}^2\\
		+(6L^4\hat{\eta}^2+3L^2\gamma^2\norm{\vI-\vW}^2)\norm{\vX^k-\vone\oX^k}^2\\
		~~+3L^2\hat{\eta}^2\norm{\vY^k-\vone\oY^k}^2+3L^2C\gamma^2\norm{\vI-\vW}^2\norm{\vX^k-\vH^k_x}^2.
	\end{multline}
	Substituting \eqref{GradplusminusGrad} and \eqref{ineqN:Yk_Ykhat} into \eqref{YkplusCGT}, we conclude that 
	\begin{align}
		\nonumber &\EE \left[\norm{\vY^{k+1} - \vone \oY^{k+1}}^2\middle|\cF^k\right] \\
		\nonumber\leq& 12nL^4\frac{\hat{\eta}^2}{s\gamma}\norm{\oX^k-\vx^*}^2
		\nonumber+(12L^4\frac{\hat{\eta}^2}{s\gamma}+\frac{6\lambda L^2\gamma}{s})\norm{\vX^k-\vone\oX^k}^2\\
		&\nonumber+(\frac{1+\tilde{\rho}^2}{2}+\frac{6L^2}{s}\frac{\hat{\eta}^2}{\gamma})\norm{\vY^k-\vone\oY^k}^2
		\\
		\label{Ykplus1}&+\frac{6C\lambda L^2\gamma}{s}\norm{\vX^k-\vH^k_x}^2
		+\frac{2C\lambda\gamma}{s}\norm{\vY^k-\vH_y^k}^2.
	\end{align}
	
	\subsubsection{Difference compression error of decision variables}
	For  convenience, denote $\cC^r=\frac{\cC}{r}$. 
	According to  the $\Call{COMM}{}$ procedure in Algorithm \ref{Alg:CGTe} for the decision variable, i.e., $\vH_x^{k+1}=\alpha_x \vH_x^k+(1-\alpha_x)\widehat{\vX}^k=\vH_x^k+\alpha_x\vQ_x^k$ and Lemma \ref{lem:UV}, we know
	\begin{align}
		\nonumber & \norm{\vX^{k+1} - \vH_x^{k+1}}^2\\
		\nonumber =& \norm{\vX^{k+1} -\vX^k+\vX^k- \vH_x^k-\alpha_x r \frac{\vQ_x^k}{r}}^2 \\
		\nonumber =& \Big\|\vX^{k+1} -\vX^k+\alpha_x r (\vX^k-\vH_x^k-\cC^r(\vX^k-\vH_x^k))\\
		\nonumber&+(1-\alpha_x r)(\vX^k- \vH_x^k)\Big\|^2\\
		\nonumber  \leq& \tau_x\Big\|\alpha_x r (\vX^k-\vH_x^k-\cC^r(\vX^k-\vH_x^k))\\
		\nonumber&+(1-\alpha_x r)(\vX^k- \vH_x^k)\Big\|^2+\frac{\tau_x}{\tau_x-1}\norm{\vX^{k+1} -\vX^k}^2\\
		\nonumber \leq& \tau_x\bigg[\alpha_x r\norm{\vX^k-\vH_x^k-\cC^r(\vX^k-\vH_x^k)}^2\\
		\label{ineq:HkxError}& +(1-\alpha_x r)\norm{\vX^k- \vH_x^k}^2\bigg]+\frac{\tau_x}{\tau_x-1}\norm{\vX^{k+1} -\vX^k}^2,
	\end{align}
	where $\alpha_x$ satisfies $0<\alpha_x r\leq 1$. The last inequality holds from the convexity of norm. 
	Taking conditional expectation on both sides of \eqref{ineq:HkxError}, we get
	\begin{align}
		\nonumber & \EE \left[\norm{\vX^{k+1} - \vH_x^{k+1}}^2\middle|\cF^k\right]\\
		\nonumber \leq& \tau_x\left[\alpha_x r(1-\delta)+(1-\alpha_x r)\right]\norm{\vX^k-\vH_x^k}^2\\
		\label{ineq:EHkxError}&+\frac{\tau_x}{\tau_x-1}\EE \left[\norm{\vX^{k+1} -\vX^k}^2\middle|\cF^k\right],
	\end{align}
	where the inequality stems from Assumption \ref{Assumption:General}. 
	For simplicity, denote $t_x=\frac{3\tau_x}{\tau_x-1}>1$ and $c_x=\tau_x\left[\alpha_x r(1-\delta)+(1-\alpha_x r)\right]=\tau_x(1-\alpha_x r \delta)<1$. Then, plugging \eqref{ineq:XkplusminusXk2} into \eqref{ineq:EHkxError}, we obtain 
	\begin{align}
		\nonumber & \EE \left[\norm{\vX^{k+1} - \vH_x^{k+1}}^2\middle|\cF^k\right]\\
		\nonumber\leq& 2t_xnL^2\hat{\eta}^2\norm{\oX^k-\vx^*}^2+(2t_xL^2\hat{\eta}^2+t_x\lambda\gamma^2)\norm{\vX^k-\vone\oX^k}^2\\
		&+t_x\hat{\eta}^2\norm{\vY^k-\vone\oY^k}^2+(c_x+t_xC\lambda\gamma^2)\norm{\vX^k-\vH^k_x}^2.
	\end{align}
	
	\subsubsection{Difference compression error of gradient trackers}
	Based on the same bounding method as above, we have
	\begin{align}
		\nonumber & \EE \left[\norm{\vY^{k+1} - \vH_y^{k+1}}^2\middle|\cF^k\right]\\
		\nonumber\leq& \tau_y\left[\alpha_y r(1-\delta)+(1-\alpha_y r)\right]\norm{\vY^k-\vH_y^k}^2\\
		\label{ineq:EHkyError}&+\frac{\tau_y}{\tau_y-1}\EE \left[\norm{\vY^{k+1} -\vY^k}^2\middle|\cF^k\right].
	\end{align}
	In light of \eqref{ineq:EHkyError}, we still need to bound
	\begin{align}
		\nonumber  &\EE \left[\norm{\vY^{k+1} -\vY^k}^2\middle|\cF^k\right]\\
		\nonumber= &\EE \bigg[\Big\|\gamma(\vI-\vW)(\vY^k-\widehat{\vY}^k)-\gamma(\vI-\vW)(\vY^k-\vone\oY^k)\\
		\nonumber&+(\nabla^{k+1}-\nabla^k)\Big\|^2\bigg|\cF^k\bigg]\\
		\nonumber\leq&
		3\gamma^2\norm{\vI-\vW}^2\EE\left[\norm{\vY^k-\widehat{\vY}^k}^2\middle|\cF^k\right]\\
		\nonumber &+3\gamma^2\norm{\vI-\vW}^2\norm{\vY^k-\vone\oY^k}^2\\
		\label{ineq:EHkyError1}&+3\EE\left[\norm{\nabla^{k+1}-\nabla^k}^2\middle|\cF^k\right].
	\end{align}
	Plugging \eqref{ineqN:Yk_Ykhat}, \eqref{GradplusminusGrad} and \eqref{ineq:EHkyError1} into \eqref{ineq:EHkyError}, we have
	\begin{align}
		\nonumber & \EE \left[\norm{\vY^{k+1} - \vH_y^{k+1}}^2\middle|\cF^k\right]\\
		\nonumber  \leq& (c_y+t_yC\lambda\gamma^2)\norm{\vY^k-\vH_y^k}^2
		+6t_ynL^4\hat{\eta}^2\norm{\oX^k-\vx^*}^2\\
		\nonumber &+(6t_yL^4\hat{\eta}^2+3t_y\lambda L^2\gamma^2)\norm{\vX^k-\vone\oX^k}^2\\
		\label{ineq:EHkyError4} \nonumber&+(t_y\lambda\gamma^2+3t_yL^2\hat{\eta}^2)\norm{\vY^k-\vone\oY^k}^2\\
		&+3t_yC\lambda L^2 \gamma^2\norm{\vX^k-\vH^k_x}^2.
	\end{align}
	where $c_y=\tau_y\left[\alpha_y r(1-\delta)+(1-\alpha_y r )\right]=\tau_y(1-\alpha_y r \delta)<1$,
	$t_y=\frac{3\tau_y}{\tau_y-1}$.
	
	\subsection{Proof of Lemma \ref{Thm:CGT}}\label{Pf:ThmCGT}
	In light of Lemma \ref{Lem:CGT}, we consider the following linear system of inequalities:
	\begin{align}\label{LIS}
		\vA\boldsymbol{\epsilon}\leq (1-\frac{1}{2}M\hat{\eta}\mu)\boldsymbol{\epsilon},
	\end{align}
	where  $\boldsymbol{\epsilon}:=[\epsilon_1,\epsilon_2,L^2\epsilon_3,\epsilon_4,L^2\epsilon_5]^{\T}$, and the elements of  $\vA$ correspond  to the coefficients in Lemma \ref{Lem:CGT}. From Lemma    \ref{implma}, if there exist an element-wise positive $\boldsymbol{\epsilon}$, then we obtain $\rho\big(\vA\big)\leq (1-\frac{1}{2}M\hat{\eta}\mu)$. 
	
	\subsubsection{First inequality in \eqref{LIS}}
	\begin{equation}\label{ineq:thm1}
		(1-\frac{3}{2}M \hat{\eta}\mu)\epsilon_1+\frac{6\hat{\eta} L^2}{\mu n}\epsilon_2+\frac{6\hat{\eta} L^2}{\mu n M}\epsilon_3 \le (1-\frac{1}{2}M\hat{\eta}\mu)\epsilon_1.
	\end{equation}
	Inequality \eqref{ineq:thm1} holds if  
	$\frac{6\hat{\eta} L^2}{\mu n}\epsilon_2 \leq \frac{6\hat{\eta} L^2}{\mu n M}\epsilon_3$ 
	and 
	$2 \frac{6\hat{\eta} L^2}{\mu n M}\epsilon_3\leq M\hat{\eta} \mu\epsilon_1$. That is,
	$\epsilon_3 \geq M \epsilon_2$ and $n\epsilon_1 \geq \frac{12\kappa^2}{M^2} \epsilon_3$, 
	where $\kappa=L/\mu$ is the condition number.
	\subsubsection{Second inequality in \eqref{LIS}}
	\begin{align}\label{ineq2:thm1}
		&\Big(\frac{1+\tilde{\rho}^2}{2}+\frac{4L^2\hat{\eta}^2}{s\gamma}\Big)\epsilon_2 
		+\frac{2\hat{\eta}^2}{s\gamma}L^2\epsilon_3
		+\frac{2C\lambda\gamma}{s}\epsilon_4
		+\frac{4nL^2\hat{\eta}^2}{s\gamma}\epsilon_1 \nonumber\\
		\leq&(1-\frac{1}{2}M\hat{\eta}\mu)\epsilon_2.
	\end{align}
	Recalling $\frac{1-\tilde{\rho}^2}{2}=(1+\tilde{\rho})\frac{1-\tilde{\rho}}{2}>\frac{1-\tilde{\rho}}{2}=\frac{\gamma s}{2}$, relation \eqref{ineq2:thm1} holds if 
	\begin{align}\label{ineq2:thm2}
		\frac{2(2n\epsilon_1+2\epsilon_2+\epsilon_3)}{s}L^2\frac{\hat{\eta}^2}{\gamma}
		+\frac{2C\lambda}{s}\gamma\epsilon_4
		\leq(\frac{s\gamma}{2}-\frac{1}{2}M\hat{\eta}\mu)\epsilon_2.
	\end{align}
	Dividing $\gamma$ on both sides of \eqref{ineq2:thm2}, we get
	\begin{align}\label{ineq2:thm3}
		\frac{M\mu\epsilon_2 }{2}\frac{\hat{\eta}}{\gamma}      +\frac{2(2n\epsilon_1+2\epsilon_2+\epsilon_3)}{s}L^2\frac{\hat{\eta}^2}{\gamma^2}
		+\frac{2C\lambda}{s}\epsilon_4
		\leq\frac{ s\epsilon_2}{2}.
	\end{align}
	It is sufficient that
	$
	3\frac{M\mu\epsilon_2 }{2}\frac{\hat{\eta}}{\gamma}\leq\frac{ s\epsilon_2}{2}$,
	$3\frac{2(2n\epsilon_1+2\epsilon_2+\epsilon_3)}{s}L^2\frac{\hat{\eta}^2}{\gamma^2}\leq\frac{ s\epsilon_2}{2}$, and
	$3\frac{2C\lambda}{s}\epsilon_4\leq\frac{ s\epsilon_2}{2}
	$.
	Therefore, if we have  $\hat{\eta}\leq \min\{\frac{s}{3M}\frac{\gamma}{\mu}, \sqrt{\frac{\epsilon_2}{12(2n\epsilon_1+2\epsilon_2+\epsilon_3)}}s\frac{\gamma}{L} \}$ and $\epsilon_2\geq \frac{12C\lambda}{s^2}\epsilon_4$, then \eqref{ineq2:thm1} can be demonstrated.
	
	\subsubsection{Third inequality in \eqref{LIS}}
	\begin{multline}\label{ineq3:thm1}
		12nL^4\frac{\hat{\eta}^2}{s\gamma}\epsilon_1
		+\Big(12L^4\frac{\hat{\eta}^2}{s\gamma}+\frac{6\lambda L^2\gamma}{s}\Big)\epsilon_2
		+\Big(\frac{1+\tilde{\rho}^2}{2}+\frac{6L^2\hat{\eta}^2}{s\gamma}\Big)L^2\epsilon_3\\
		+\frac{6C\lambda L^2 \gamma }{s}\epsilon_4
		+\frac{2C\lambda \gamma }{s}L^2\epsilon_5
		\leq(1-\frac{1}{2}M\hat{\eta}\mu)L^2\epsilon_3.
	\end{multline}
	Dividing $L^2$ on the both side of \eqref{ineq3:thm1}, we have
	\begin{multline}\label{ineq3:thm01}
		12nL^2\frac{\hat{\eta}^2}{s\gamma}\epsilon_1
		+\Big(12L^2\frac{\hat{\eta}^2}{s\gamma}+\frac{6\lambda \gamma}{s}\Big)\epsilon_2
		+\Big(\frac{1+\tilde{\rho}^2}{2}+\frac{6L^2\hat{\eta}^2}{s\gamma}\Big)\epsilon_3\\
		+\frac{6C\lambda  \gamma }{s}\epsilon_4
		+\frac{2C\lambda \gamma }{s}\epsilon_5
		\leq(1-\frac{1}{2}M\hat{\eta}\mu)\epsilon_3.
	\end{multline}
	Based on arguments similar those for deriving the second inequality, the third inequality holds if 
	\begin{align}\label{ineq3:thm2}
		&\frac{M\mu\epsilon_3}{2}\frac{\hat{\eta}}{\gamma}+\frac{6(2n\epsilon_1+2\epsilon_2+\epsilon_3)}{s}L^2\frac{\hat{\eta}^2}{\gamma^2} \nonumber\\
		&+\frac{2\lambda\left(3\epsilon_2+C(3\epsilon_4+\epsilon_5)\right)}{s} \leq \frac{s\epsilon_3}{2}.
	\end{align}
	It is sufficient that
	$
	3M\mu\epsilon_3\frac{\hat{\eta}}{\gamma}\leq\frac{s\epsilon_3}{2}$,
	$3\frac{6(2n\epsilon_1+2\epsilon_2+\epsilon_3)}{s}L^2\frac{\hat{\eta}^2}{\gamma^2}\leq\frac{s\epsilon_3}{2}$,	
	and $3\frac{2\lambda\left(3\epsilon_2+C(3\epsilon_4+\epsilon_5)\right)}{s}\leq\frac{s\epsilon_3}{2}
	$.
	Thus, if there holds $\epsilon_3\geq\frac{12\lambda\left(3\epsilon_2+C(3\epsilon_4+\epsilon_5)\right)}{s^2}$ and $\hat{\eta}\leq \min\{ \frac{s}{3M}\frac{\gamma}{\mu}, \sqrt{\frac{\epsilon_3}{36(2n\epsilon_1+2\epsilon_2+\epsilon_3)}}s\frac{\gamma}{L}\}$, we can demonstrate \eqref{ineq3:thm1}. 					
	\subsubsection{Fourth inequality in \eqref{LIS}}
	\begin{multline} \label{ineq4:thm1}
		2nt_xL^2\hat{\eta}^2\epsilon_1
		+(2t_xL^2\hat{\eta}^2
		+t_x\lambda\gamma^2)\epsilon_2\\
		+t_x\hat{\eta}^2 L^2\epsilon_3
		+(c_x+t_xC\lambda\gamma^2)\epsilon_4
		\leq(1-\frac{1}{2}M\hat{\eta}\mu)\epsilon_4.
	\end{multline}
	It is equivalent to
	\begin{multline} \label{ineq4:thm2}
		t_x(2n\epsilon_1+2\epsilon_2+\epsilon_3)L^2\hat{\eta}^2\\
		+(t_x\lambda\epsilon_2+t_xC\lambda\epsilon_4)\gamma^2
		+\frac{1}{2}M\hat{\eta}\mu\epsilon_4
		\leq(1-c_x)\epsilon_4.
	\end{multline}
	Inequality \eqref{ineq4:thm2} holds if $\hat{\eta}\leq\frac{\gamma}{L}$ and $\gamma\leq\min\{1,\frac{1-c_x}{m_x}\epsilon_4\}$,
	where 
	$
	m_x:=t_x(2n\epsilon_1+2\epsilon_2+\epsilon_3)+t_x\lambda(\epsilon_2+C\epsilon_4)
	+\frac{M\epsilon_4}{2\kappa}.
	$
	\subsubsection{Fifth inequality in \eqref{LIS}}
	\begin{multline}\label{ineq5:thm1}
		6nt_yL^4\hat{\eta}^2\epsilon_1
		+(6t_yL^4\hat{\eta}^2+3t_y\lambda L^2\gamma^2)\epsilon_2
		+(t_y\lambda \gamma^2+3t_yL^2\hat{\eta}^2)L^2\epsilon_3\\
		+3t_yC\lambda L^2\gamma^2\epsilon_4
		+(c_y+t_yC\lambda \gamma^2)L^2\epsilon_5
		\leq(1-\frac{1}{2}M\hat{\eta}\mu)L^2\epsilon_5.
	\end{multline}
	Dividing $L^2$ on both sides of \eqref{ineq5:thm1}, we have
	\begin{multline}\label{ineq5:thm2}
		3t_y(2n\epsilon_1+2\epsilon_2+\epsilon_3)L^2\hat{\eta}^2
		+t_y\lambda(3\epsilon_2+\epsilon_3\\
		\qquad+3C\epsilon_4+C\epsilon_5)\gamma^2
		+\frac{M\mu\epsilon_5}{2}\hat{\eta} 
		\leq(1-c_y)\epsilon_5.
	\end{multline}
	Inequality \eqref{ineq5:thm2}  holds if $\hat{\eta}\leq\frac{\gamma}{L}$ and $\gamma\leq\min\{1,\frac{1-c_y}{m_y}\epsilon_5\}$,
	where 
	$
	m_y:=3t_y(2n\epsilon_1+2\epsilon_2+\epsilon_3)
	+t_y\lambda(3\epsilon_2+\epsilon_3
	+3C\epsilon_4+C\epsilon_5)
	+\frac{M\epsilon_5}{2\kappa}.
	$
	In short,  if the positive constants $\epsilon_1$-$\epsilon_5$, consensus step-size $\gamma$ and step-size $\hat{\eta}$ satisfy the following conditions,	
	\begin{align}
		n\epsilon_1 \geq  &\frac{12\kappa^2}{M^2}\epsilon_3,
		\epsilon_2\geq \frac{12C\lambda}{s^2}\epsilon_4,
		\epsilon_3\geq M\epsilon_2, \nonumber \\
		\epsilon_3\geq &\frac{12\lambda\left(3\epsilon_2+C(3\epsilon_4+\epsilon_5)\right)}{s^2},
		\epsilon_4>0,
		\epsilon_5>0,
		\label{condition_epsilon}\\
		\gamma\leq &\min\Big\{1,
		\frac{1-c_x}{m_x}\epsilon_4,
		\frac{1-c_y}{m_y}\epsilon_5\Big\},\label{gammaCGT}\\
		\hat{\eta}\le &\min\Big\{ \frac{s}{3M}\frac{\gamma}{\mu},
		\frac{\gamma}{L},
		\sqrt{\frac{\epsilon_2}{12(2n\epsilon_1+2\epsilon_2+\epsilon_3)}}\frac{s\gamma}{L},\nonumber\\
		&\qquad\qquad\sqrt{\frac{\epsilon_3}{36(2n\epsilon_1+2\epsilon_2+\epsilon_3)}}\frac{s\gamma}{L},
		\Big\},\label{etaCGT}
	\end{align}
	we can establish the linear system of inequalities in (\ref{LIS}).  
	For \eqref{condition_epsilon}, it is easy to verify that there exist solutions to $\epsilon_1$-$\epsilon_5$. Noticing that 	     $\sqrt{\frac{\epsilon_2}{12(2n\epsilon_1+2\epsilon_2+\epsilon_3)}}\frac{s\gamma}{L}
	=\sqrt{\frac{1}{12(2n\epsilon_1/\epsilon_2+2+\epsilon_3/\epsilon_2)}}\frac{s\gamma}{L}$ 
	and 
	$\sqrt{\frac{\epsilon_3}{36(2n\epsilon_1+2\epsilon_2+\epsilon_3)}}\frac{s\gamma}{L}
	=\sqrt{\frac{1}{36(2n\epsilon_1/\epsilon_3+2\epsilon_2/\epsilon_3+1)}}\frac{s\gamma}{L}$
	are both less than $\frac{s}{3M}\frac{\gamma}{\mu}$ and $\frac{\gamma}{L}$,  we obtain the upper bound on the maximum step-size.

	\subsection{Proof of Theorem \ref{Thm1:C-GT}}\label{Pf:Thm1CGT}
	Taking $\tau_x=\tau_y=\frac{1}{1-0.5\delta}$ and noticing that $\alpha_x=\alpha_y=1/r$, we have
	$t_x=t_y=\frac{6}{\delta}$ and $c_x=c_y=\frac{1-\delta}{1-0.5\delta}$. 
	Recalling the  relations of $\epsilon_1$-$\epsilon_5$ in Lemma \ref{Thm:CGT}, we can take $\epsilon_4=\epsilon_5=1$, $\epsilon_2=\frac{12C\lambda}{s^2}$, $\epsilon_3=\frac{432C\lambda^2}{s^4}+\frac{48C\lambda}{s^2}\geq4\epsilon_2$ and $n\epsilon_1=\frac{12\kappa^2}{M^2}\Big(\frac{432C\lambda^2}{s^4}+\frac{48C\lambda}{s^2}\Big)$. Then, we know 
	$\frac{\epsilon_2}{12(2n\epsilon_1+2\epsilon_2+\epsilon_3)}\leq \frac{\epsilon_3}{36(2n\epsilon_1+2\epsilon_2+\epsilon_3)}$.
	Meanwhile, noticing that 
	$
	m_y:=3\frac{6}{\delta}[(2n\epsilon_1+2\epsilon_2+\epsilon_3)
	+\lambda(\epsilon_2+C)]+\frac{M}{2\kappa}+\frac{6}{\delta}\lambda(\epsilon_3+C)
	$ 
	and 
	$ m_x=\frac{6}{\delta}[(2n\epsilon_1+2\epsilon_2+\epsilon_3)
	+\lambda(\epsilon_2+C)]+\frac{M}{2\kappa}
	$, we get $m_y\geq m_x$. For simplicity, denote $m=m_y$ where the specific values for the constants $\epsilon_1$-$\epsilon_5$ are given above. Then from Lemma \ref{Thm:CGT}, we obtain the upper bounds on the consensus step-size and the maximum step-size, which completes the proof.

\end{document}

%% file: notation.tex
\newcommand{\ve}{{\mathbf{e}}}

\newcommand{\vv}{{\mathbf{v}}}
\newcommand{\vw}{{\mathbf{w}}}
\newcommand{\vx}{{\mathbf{x}}}
\newcommand{\vy}{{\mathbf{y}}}

\newcommand{\vA}{{\mathbf{A}}}

\newcommand{\vD}{{\mathbf{D}}}
\newcommand{\vE}{{\mathbf{E}}}
\newcommand{\vF}{{\mathbf{F}}}

\newcommand{\vH}{{\mathbf{H}}}
\newcommand{\vI}{{\mathbf{I}}}

\newcommand{\vM}{{\mathbf{M}}}

\newcommand{\vQ}{{\mathbf{Q}}}

\newcommand{\vU}{{\mathbf{U}}}
\newcommand{\vV}{{\mathbf{V}}}
\newcommand{\vW}{{\mathbf{W}}}
\newcommand{\vX}{{\mathbf{X}}}
\newcommand{\vY}{{\mathbf{Y}}}
\newcommand{\vZ}{{\mathbf{Z}}}

\newcommand{\cC}{{\mathcal{C}}}

\newcommand{\cF}{{\mathcal{F}}}

\newcommand{\cO}{{\mathcal{O}}}

\newcommand{\cQ}{{\mathcal{Q}}}

\newcommand{\EE}{{\mathbb{E}}}
\newcommand{\RR}{\mathbb{R}}
\newcommand{\vzero}{\mathbf{0}}
\newcommand{\vone}{{\mathbf{1}}}
\newcommand{\oX}{\overline{\vX}}
\newcommand{\oY}{\overline{\vY}}
\newcommand{\T}{\intercal}
\newcommand{\tW}{\widetilde{\vW}}
\newcommand{\norm}[1]{\left\| #1 \right\|}      % norm 

%% file: main_arxiv.bbl
% Generated by IEEEtran.bst, version: 1.14 (2015/08/26)
\begin{thebibliography}{10}
\providecommand{\url}[1]{#1}
\csname url@samestyle\endcsname
\providecommand{\newblock}{\relax}
\providecommand{\bibinfo}[2]{#2}
\providecommand{\BIBentrySTDinterwordspacing}{\spaceskip=0pt\relax}
\providecommand{\BIBentryALTinterwordstretchfactor}{4}
\providecommand{\BIBentryALTinterwordspacing}{\spaceskip=\fontdimen2\font plus
\BIBentryALTinterwordstretchfactor\fontdimen3\font minus
  \fontdimen4\font\relax}
\providecommand{\BIBforeignlanguage}[2]{{%
\expandafter\ifx\csname l@#1\endcsname\relax
\typeout{** WARNING: IEEEtran.bst: No hyphenation pattern has been}%
\typeout{** loaded for the language `#1'. Using the pattern for}%
\typeout{** the default language instead.}%
\else
\language=\csname l@#1\endcsname
\fi
#2}}
\providecommand{\BIBdecl}{\relax}
\BIBdecl

\bibitem{Cohen2016Distributed}
K.~Cohen, A.~Nedi{\'c}, and R.~Srikant, ``Distributed learning algorithms for
  spectrum sharing in spatial random access wireless networks,'' \emph{IEEE
  Transactions on Automatic Control}, vol.~62, no.~6, pp. 2854--2869, 2016.

\bibitem{Nedic2018Distributed}
A.~Nedi{\'c} and J.~Liu, ``Distributed optimization for control,'' \emph{Annual
  Review of Control, Robotics, and Autonomous Systems}, vol.~1, no.~1, pp.
  77--103, 2018.

\bibitem{Nedic2020Distributed}
A.~Nedi{\'c}, ``Distributed gradient methods for convex machine learning
  problems in networks: {{Distributed}} optimization,'' \emph{IEEE Signal
  Processing Magazine}, vol.~37, no.~3, pp. 92--101, 2020.

\bibitem{Nedic2009distributed}
A.~Nedi{\'c} and A.~Ozdaglar, ``Distributed subgradient methods for multi-agent
  optimization,'' \emph{IEEE Transactions on Automatic Control}, vol.~54,
  no.~1, pp. 48--61, Jan. 2009.

\bibitem{Shi2015Extra}
W.~Shi, Q.~Ling, G.~Wu, and W.~Yin, ``{EXTRA}: An exact first-order algorithm
  for decentralized consensus optimization,'' \emph{SIAM Journal on
  Optimization}, vol.~25, no.~2, pp. 944--966, 2015.

\bibitem{Xu2015Augmented}
J.~Xu, S.~Zhu, Y.~C. Soh, and L.~Xie, ``Augmented distributed gradient methods
  for multi-agent optimization under uncoordinated constant stepsizes,'' in
  \emph{Proceedings of the 54th {IEEE} Conference on Decision and Control
  ({CDC})}.\hskip 1em plus 0.5em minus 0.4em\relax {IEEE}, 2015, pp.
  2055--2060.

\bibitem{Di2016Next}
P.~Di~Lorenzo and G.~Scutari, ``{NEXT}: {In}-network nonconvex optimization,''
  \emph{IEEE Transactions on Signal and Information Processing over Networks},
  vol.~2, no.~2, pp. 120--136, 2016.

\bibitem{Nedic2017achieving}
A.~Nedi{\'c}, A.~Olshevsky, and W.~Shi, ``Achieving geometric convergence for
  distributed optimization over time-varying graphs,'' \emph{SIAM Journal on
  Optimization}, vol.~27, no.~4, pp. 2597--2633, 2017.

\bibitem{Qu2018Harnessing}
G.~Qu and N.~Li, ``Harnessing smoothness to accelerate distributed
  optimization,'' \emph{IEEE Transactions on Control of Network Systems},
  vol.~5, no.~3, pp. 1245--1260, Sep. 2018.

\bibitem{Nedic2017Geometrically}
A.~Nedi{\'c}, A.~Olshevsky, W.~Shi, and C.~A. Uribe, ``Geometrically convergent
  distributed optimization with uncoordinated step-sizes,'' in \emph{2017
  American Control Conference ({ACC})}.\hskip 1em plus 0.5em minus 0.4em\relax
  {IEEE}, 2017, pp. 3950--3955.

\bibitem{Pu2020distributed}
S.~Pu and A.~Nedi{\'c}, ``Distributed stochastic gradient tracking methods,''
  \emph{Mathematical Programming}, pp. 1--49, 2020.

\bibitem{Li2020communication}
B.~Li, S.~Cen, Y.~Chen, and Y.~Chi, ``Communication-efficient distributed
  optimization in networks with gradient tracking and variance reduction,'' in
  \emph{International Conference on Artificial Intelligence and
  Statistics}.\hskip 1em plus 0.5em minus 0.4em\relax PMLR, 2020, pp.
  1662--1672.

\bibitem{Nedic2015Distributed}
A.~Nedi{\'c} and A.~Olshevsky, ``Distributed optimization over time-varying
  directed graphs,'' \emph{IEEE Transactions on Automatic Control}, vol.~60,
  no.~3, pp. 601--615, 2015.

\bibitem{Xie2018Distributed}
P.~Xie, K.~You, R.~Tempo, S.~Song, and C.~Wu, ``Distributed convex optimization
  with inequality constraints over time-varying unbalanced digraphs,''
  \emph{IEEE Transactions on Automatic Control}, vol.~63, no.~12, pp.
  4331--4337, 2018.

\bibitem{sun2022distributed}
Y.~Sun, G.~Scutari, and A.~Daneshmand, ``Distributed optimization based on
  gradient tracking revisited: Enhancing convergence rate via surrogation,''
  \emph{SIAM Journal on Optimization}, vol.~32, no.~2, pp. 354--385, 2022.

\bibitem{Tsianos2012PushSum}
K.~I. Tsianos, S.~Lawlor, and M.~G. Rabbat, ``Push-sum distributed dual
  averaging for convex optimization,'' in \emph{Proceedings of the 51st {IEEE}
  Conference on Decision and Control ({CDC})}.\hskip 1em plus 0.5em minus
  0.4em\relax {IEEE}, 2012, pp. 5453--5458.

\bibitem{Xin2018Linear}
R.~Xin and U.~A. Khan, ``A linear algorithm for optimization over directed
  graphs with geometric convergence,'' \emph{IEEE Control Systems Letters},
  vol.~2, no.~3, pp. 315--320, 2018.

\bibitem{Pu2021Push}
S.~Pu, W.~Shi, J.~Xu, and A.~Nedi{\'c}, ``Push\textendash{Pull} gradient
  methods for distributed optimization in networks,'' \emph{IEEE Transactions
  on Automatic Control}, vol.~66, no.~1, pp. 1--16, Jan. 2021.

\bibitem{Xin2020General}
R.~Xin, S.~Pu, A.~Nedi{\'c}, and U.~A. Khan, ``A general framework for
  decentralized optimization with first-order methods,'' \emph{Proceedings of
  the IEEE}, vol. 108, no.~11, pp. 1869--1889, 2020.

\bibitem{Pu2020Robust}
S.~Pu, ``A robust gradient tracking method for distributed optimization over
  directed networks,'' in \emph{Proceedings of the 59th {IEEE} Conference on
  Decision and Control ({CDC})}.\hskip 1em plus 0.5em minus 0.4em\relax {IEEE},
  2020, pp. 2335--2341.

\bibitem{1-bit-sgd}
F.~Seide, H.~Fu, J.~Droppo, G.~Li, and D.~Yu, ``1-bit stochastic gradient
  descent and application to data-parallel distributed training of speech
  {{DNNs}},'' in \emph{Interspeech 2014}, Sep. 2014.

\bibitem{NIPS2017_qsgd}
D.~Alistarh, D.~Grubic, J.~Li, R.~Tomioka, and M.~Vojnovic, ``{QSGD}:
  Communication-efficient {SGD} via gradient quantization and encoding,'' in
  \emph{Advances in Neural Information Processing Systems}, 2017, pp.
  1709--1720.

\bibitem{signSGD-ICML18}
J.~Bernstein, Y.-X. Wang, K.~Azizzadenesheli, and A.~Anandkumar, ``{SIGNSGD}:
  Compressed optimisation for non-convex problems,'' in \emph{Proceedings of
  the 35th International Conference on Machine Learning}, 2018, pp. 559--568.

\bibitem{Stich2018Sparsified}
S.~U. Stich, J.-B. Cordonnier, and M.~Jaggi, ``Sparsified {SGD} with memory,''
  in \emph{Advances in Neural Information Processing Systems}, 2018, pp.
  4452--4463.

\bibitem{karimireddy2019Error}
S.~P. Karimireddy, Q.~Rebjock, S.~Stich, and M.~Jaggi, ``Error feedback fixes
  signsgd and other gradient compression schemes,'' in \emph{Proceedings of the
  36th International Conference on Machine Learning}.\hskip 1em plus 0.5em
  minus 0.4em\relax {PMLR}, 2019, pp. 3252--3261.

\bibitem{mishchenko2019distributed}
K.~Mishchenko, E.~Gorbunov, M.~Tak{\'a}{\v c}, and P.~Richt{\'a}rik,
  ``Distributed learning with compressed gradient differences,'' \emph{arXiv
  preprint arXiv:1901.09269}, 2019.

\bibitem{tang2019doublesqueeze}
H.~Tang, C.~Yu, X.~Lian, T.~Zhang, and J.~Liu, ``{DoubleSqueeze}: {Parallel}
  stochastic gradient descent with double-pass error-compensated compression,''
  in \emph{Proceedings of the 36th International Conference on Machine
  Learning}, 2019, pp. 6155--6165.

\bibitem{Stich2020Communication}
S.~U. Stich, ``On communication compression for distributed optimization on
  heterogeneous data,'' \emph{arXiv preprint arXiv:2009.02388}, 2020.

\bibitem{Beznosikov2020Biased}
A.~Beznosikov, S.~Horv{\'a}th, P.~Richt{\'a}rik, and M.~Safaryan, ``On biased
  compression for distributed learning,'' \emph{arXiv preprint
  arXiv:2002.12410}, 2020.

\bibitem{Xu2020Compressed}
H.~Xu, C.-Y. Ho, A.~M. Abdelmoniem, A.~Dutta, E.~H. Bergou, K.~Karatsenidis,
  M.~Canini, and P.~Kalnis, ``Compressed communication for distributed deep
  learning: {Survey} and quantitative evaluation,'' {KAUST}, Tech. Rep., 2020.

\bibitem{tang_NIPS2018_7992}
H.~Tang, S.~Gan, C.~Zhang, T.~Zhang, and J.~Liu, ``Communication compression
  for decentralized training,'' in \emph{Advances in Neural Information
  Processing Systems}, 2018, pp. 7652--7662.

\bibitem{ksj2019choco}
A.~Koloskova, S.~U. Stich, and M.~Jaggi, ``Decentralized stochastic
  optimization and gossip algorithms with compressed communication,'' in
  \emph{Proceedings of the 36th International Conference on Machine
  Learning}.\hskip 1em plus 0.5em minus 0.4em\relax {PMLR}, 2019, pp.
  3479--3487.

\bibitem{koloskova*2020decentralized}
A.~Koloskova, T.~Lin, S.~U. Stich, and M.~Jaggi, ``Decentralized deep learning
  with arbitrary communication compression,'' in \emph{International Conference
  on Learning Representations}, 2020.

\bibitem{liu2020linear}
X.~Liu, Y.~Li, R.~Wang, J.~Tang, and M.~Yan, ``Linear convergent decentralized
  optimization with compression,'' in \emph{International Conference on
  Learning Representations}, 2020.

\bibitem{li2019decentralized}
Z.~Li, W.~Shi, and M.~Yan, ``A decentralized proximal-gradient method with
  network independent step-sizes and separated convergence rates,'' \emph{IEEE
  Transactions on Signal Processing}, vol.~67, no.~17, pp. 4494--4506, 2019.

\bibitem{kovalev2021linearly}
D.~Kovalev, A.~Koloskova, M.~Jaggi, P.~Richtarik, and S.~Stich, ``A linearly
  convergent algorithm for decentralized optimization: Sending less bits for
  free!'' in \emph{International Conference on Artificial Intelligence and
  Statistics}.\hskip 1em plus 0.5em minus 0.4em\relax PMLR, 2021, pp.
  4087--4095.

\bibitem{Lee2021finite}
C.-S. Lee, N.~Michelusi, and G.~Scutari, ``Finite-bit quantization for
  distributed algorithms with linear convergence,'' \emph{arXiv preprint
  arXiv:2107.11304}, 2021.

\bibitem{Kajiyama2020Linear}
Y.~Kajiyama, N.~Hayashi, and S.~Takai, ``Linear convergence of consensus-based
  quantized optimization for smooth and strongly convex cost functions,''
  \emph{IEEE Transactions on Automatic Control}, vol.~66, no.~3, pp.
  1254--1261, 2021.

\bibitem{Basu2020qsparse}
D.~Basu, D.~Data, C.~Karakus, and S.~N. Diggavi, ``{Qsparse-Local-SGD}:
  Distributed sgd with quantization, sparsification, and local computations,''
  \emph{IEEE Journal on Selected Areas in Information Theory}, vol.~1, no.~1,
  pp. 217--226, 2020.

\bibitem{Horn2012Matrix}
R.~A. Horn and C.~R. Johnson, \emph{Matrix Analysis}, 2nd~ed.\hskip 1em plus
  0.5em minus 0.4em\relax {Cambridge university press}, 2012.

\end{thebibliography}
